\documentclass[final,onefignum,onetabnum]{siamart220329}

\usepackage{lipsum}
\usepackage{amsfonts}
\usepackage{graphicx}
\usepackage{epstopdf}
\usepackage{algorithmic}
\usepackage{amsmath}
\usepackage{dutchcal} 
\usepackage{amsopn}
\usepackage{verbatim} 
\usepackage{enumitem}

\ifpdf
\DeclareGraphicsExtensions{.eps,.pdf,.png,.jpg}
\else
\DeclareGraphicsExtensions{.eps}
\fi

\usepackage{xcolor}
\definecolor{darkgreen}{rgb}{0.1,0.6,0.1}

\newcommand{\vv}{\mathbf{v}}  
\newcommand{\ee}{E} 
\renewcommand{\S}{S} 
\newcommand{\J}{\mathbf{J}}  
\newcommand{\A}{\mathbf{A}}  
\newcommand{\bs}{\boldsymbol{\sigma}}  
\newcommand{\bbeta}{\boldsymbol{\beta}}  
\newcommand{\bG}{\boldsymbol{\Gamma}}  
\newcommand{\mde}[2]{\frac {d #1}{d #2}}
\newcommand{\halb}{\frac{1}{2}}
\newcommand{\x}{\mathbf{x}}  
\newcommand{\q}{\mathbf{q}} 
\newcommand{\w}{\mathbf{w}} 
\newcommand{\n}{\mathbf{n}} 

\newcommand{\M}{\mathbf{M}}
\newcommand{\tr}{\textnormal{tr}}
\newcommand{\lnpc}{l_{pc} \n_{pc}} 
\newcommand{\Lx}{\boldsymbol{\xi}}
\newcommand{\detA}{\text{det}(\A)}

\newcommand{\auxf}{\lambda}  

\allowdisplaybreaks

\newsiamremark{remark}{Remark}
\newsiamremark{hypothesis}{Hypothesis}
\crefname{hypothesis}{Hypothesis}{Hypotheses}
\newsiamthm{claim}{Claim}

\headers{Structure-preserving HTC scheme for continuum mechanics}{W. Boscheri, M. Dumbser, R. Loub\`ere and P.H. Maire}    %

\title{A structure-preserving and thermodynamically compatible cell-centered Lagrangian finite volume scheme for continuum mechanics}

\author{
	Walter Boscheri \thanks{Laboratoire de Mathématiques UMR 5127 CNRS, Universit{\'e} Savoie Mont Blanc, 73376 Le Bourget du Lac, France (\email{walter.boscheri@univ-smb.fr}) and Department of Mathematics and Computer Science, University of Ferrara, via Machiavelli 30, 44121 Ferrara, Italy}
	\and 	
	Michael Dumbser\thanks{Department of Civil, Environmental and Mechanical Engineering, University of Trento, Via Mesiano 77, 38123 Trento, Italy (\email{michael.dumbser@unitn.it})}
	\and
	Rapha\"el Loub\`ere \thanks{Institut de Math\'ematiques de Bordeaux (IMB), UMR 5251, Universit\'e de Bordeaux, CNRS, Bordeaux INP, Talence, F-33400, France (\email{raphael.loubere@u-bordeaux.fr})}	 	
	\and
	Pierre-Henri Maire \thanks{CEA CESTA, 33116 Le Barp, France (\email{pierre-henri.maire@cea.fr})} 
}

\begin{document}

\maketitle

\begin{abstract}
In this work we present a novel structure-preserving scheme for the compatible discretization of the  Godunov-Peshkov-Romenski (GPR) model of continuum mechanics written in Lagrangian form. The governing equations fall into the larger class of overdetermined hyperbolic and thermodynamically compatible (HTC) systems of partial differential equations (PDE). This model admits an extra conservation law for the total energy (first principle of thermodynamics) and satisfies the entropy inequality (second principle of thermodynamics). Furthermore, in the absence of algebraic source terms, the distortion field of the continuum and the specific thermal impulse satisfy a curl-free condition, provided the initial data are curl-free. Last but not least, the determinant of the distortion field is related to the density of the medium, i.e. the system is also endowed with a nonlinear algebraic constraint. In the stiff relaxation limit, the system tends to the compressible Navier-Stokes equations, i.e. the GPR model is able to describe at the same time the dynamics of nonlinear solids as well as the one of fluids.      
The objective of this work is to construct and analyze a new semi-discrete thermodynamically compatible cell-centered Lagrangian finite volume scheme on moving unstructured meshes that satisfies the following structural properties of the governing PDE exactly at the discrete level: i) compatibility with the first law of thermodynamics, i.e. discrete total energy conservation; ii) compatibility with the second law of thermodynamics, i.e. discrete entropy inequality; iii) exact discrete compatibility between the density and the determinant of the distortion field; iv) exact preservation of the curl-free property of the distortion field and of the specific thermal impulse in the absence of algebraic source terms. We will show that it is indeed possible to achieve all above properties simultaneously. Unlike in existing schemes, we choose to directly discretize the entropy inequality, hence obtaining total energy conservation as a consequence of an appropriate and thermodynamically compatible discretization of all the other equations. From this discretization, property ii) above follows trivially by construction, while i) leads to provable nonlinear stability, which is an important feature, in particular for the complex PDE system under consideration here. 

The thermodynamic compatibility and thus nonlinear stability in the sense of total energy conservation is achieved via a very simple and general approach recently introduced by Abgrall \textit{et al.} by using a scalar correction factor that is defined at the nodes of the grid. This perfectly fits into the formalism of nodal solvers which is typically adopted in cell-centered Lagrangian finite volume methods. The new scheme is run on some academic benchmark problems for computational fluid and solid mechanics to show that the properties also hold in the practical implementation of the scheme. 

\end{abstract}



\section{Introduction}\label{sec.intro}

In \cite{God1961} Godunov discovered for the first time the connection between symmetric hyperbolicity in the sense of Friedrichs \cite{FriedrichsSymm} and thermodynamic compatibility with the first and second laws of thermodynamics. 10 years later, the same connection was rediscovered by Friedrichs \& Lax in  \cite{FriedrichsLax}. Subsequently, Godunov and collaborators extended this discovery to the more general formalism of symmetric hyperbolic and thermodynamically compatible (SHTC) systems, including compressible magnetohydrodynamics (MHD), nonlinear hyperelasticity and compressible multi-phase flows, see e.g. \cite{God1972MHD,GodRom1972,GodRom2003,Rom1998,RomenskiTwoPhase2010}. 
In his seminal papers \cite{Tadmor1,Tadmor2003} Tadmor was the first to achieve thermodynamic compatibility also on the discrete level inside a numerical scheme, obtaining the discrete compatibility with the \textit{entropy inequality} (second principle of thermodynamics) as a consequence of the compatible numerical scheme, while directly discretizing the \textit{total energy} conservation law (first principle of thermodynamics). Further ongoing research on this subject can be found, for example, in \cite{Tadmor2003,Hiltebrand2014,ray2016,FjordholmMishraTadmor,Ranocha2020,Chatterjee2020,ShuEntropyMHD1,ShuEntropyMHD2,GassnerEntropyGLM}. Entropy-compatible discretizations for non-conservative hyperbolic equations were forwarded, for instance, in \cite{Fjordholm2012,AbgrallBT2018,Abgrall2018}. 
For hydrodynamics and continuum mechanics a first attempt to achieve discrete total energy conservation as a consequence of all other equations rather than a discrete entropy inequality was made in \cite{HTCGPR}.   
In the context of updated Lagrangian schemes for hydrodynamics on moving meshes, total energy conservation can be obtained as a consequence of the discrete mass, momentum and internal energy equations, see e.g. \cite{CaramanaCompatible1,CaramanaCompatible2}. Further developments on total energy conserving and entropy-compatible Lagrangian schemes for hydrodynamics can be found in \cite{Despres2009,EUCCLHYD,BRAEUNIG2016127,Maire2020,HTCLagrange}. 

From a mathematical point of view, a particularly interesting system is the system of nonlinear hyperelasticity, which was written in Eulerian coordinates in \cite{GodRom1972}. It was later extended to incorporate also strain relaxation source terms that allow to describe visco-plastic solids and viscous fluids in \cite{PeshRom2014,GPRmodel}. Compared to standard hydrodynamics, it includes two additional quantities that need to be evolved in time, namely the distorsion field $\A$, which in absence of source terms is the Eulerian gradient of the Lagrangian reference coordinates, and the specific thermal impulse $\J$, which allows for hyperbolic heat conduction. 
The PDE system is symmetrizable \cite{Rom1998}, but the energy potential is not strictly convex, hence the Hessian of the energy potential is obviously symmetric, but not strictly positive definite everywhere.
However, the system has the following interesting additional mathematical properties: 
\begin{enumerate} 
	\item the system is thermodynamically compatible, i.e. it satisfies total energy conservation and 
	\item it satisfies an entropy inequality;
	\item the system admits an algebraic constraint that connects the determinant of the distorsion field $\A$ with the mass density $\rho$ of the medium; 
	\item in the absence of source terms, the distorsion field $\A$ and the thermal impulse $\J$ satisfy the stationary differential constraints (involutions) $\nabla \times \A = 0$ and $\nabla \times \J = 0$.
\end{enumerate}

The main objective of this paper is to construct a numerical scheme that satisfies simultaneously all the properties 1)-4) exactly also on the discrete level. To the best of our knowledge, up to now there exists no such scheme. Previous numerical methods were either only compatible with the involutions \cite{SIGPR}, or only thermodynamically compatible \cite{HTCGPR,HTCAbgrall}, or thermodynamically compatible and compatible with the determinant constraint \cite{BoscheriGPRGCL}, but none was able to satisfy all properties 1)-4) simultaneously. A very special feature of our scheme is that we discretize the entropy inequality directly and obtain total energy conservation as a consequence, similar to what was done in \cite{HTCGPR,HTCAbgrall,HTCLagrange}. 
The rest of this paper is organized as follows. In Section \ref{sec.model} we introduce the governing partial differential equations. In Section \ref{sec.method} we present our new thermodynamically compatible Lagrangian HTC scheme and prove that all mathematical properties 1)-4) also hold exactly on the semi-discrete level. Numerical results for some elementary test cases are shown in Section \ref{sec.results}. The paper closes with some conclusions and an outlook to future research in Section \ref{sec.concl}.

\section{A unified first order thermodynamically compatible model of continuum mechanics} \label{sec.model}
Using the material derivative $\mde{}{t} = \frac{\partial }{\partial t} + \mathbf{v} \cdot \nabla $ the unified first order hyperbolic model  of continuum mechanics of Godunov, Peshkov and Romenski (GPR)  \cite{GodRom1972,Rom1998,PeshRom2014,GPRmodel,HTCGPR}, that is able to describe the dynamics of nonlinear elastic solids at large deformations and viscous fluids with heat conduction, can be written in the updated Lagrangian form \cite{Boscheri_hyperelast_22,LGPR} as follows:
\begin{subequations}
	\label{eqn.PDE_entropy}
	\begin{align}
		&\rho \mde{\tau}{t} - \frac{\partial  v_k}{\partial x_k} =0,    \label{eqn.cl_tau}\\
		&\rho \mde{v_i}{t}  + \frac{\partial}{\partial x_k} \left( p \delta_{ik} + \sigma_{ik} \right) = 0,    \label{eqn.cl_v}\\
		&\rho \mde{S}{t}    + \frac{\partial (\rho \beta_k)}{\partial x_k}  =  
		\frac{\Gamma_{ij} \Gamma_{ij}}{T \theta_1(\tau_1)}  +
		\frac{\beta_{i}  \beta_{i}}{T \theta_2(\tau_2)}  
		\, \geq 0, \label{eqn.cl_S} \\
		&\rho \mde{A_{ik}}{t}   + \rho A_{il} \frac{\partial v_l}{\partial x_k}  = -\frac{1}{\theta_1(\tau_1)} \Gamma_{ij}, \label{eqn.cl_A}	\\ 
		&\rho \mde{J_i}{t}  + \rho J_k \frac{\partial v_k}{\partial x_i} + \rho \frac{\partial T}{\partial x_k} = -\frac{1}{\theta_2(\tau_2)} \beta_{i}. \label{eqn.cl_J} 
	\end{align}
\end{subequations}
The above model is supplemented by the following extra conservation law for the total energy, see \cite{Rom1998,PeshRom2014,GPRmodel,HTCGPR}: 
\begin{equation}
	\rho \mde{E}{t} + \frac{\partial}{\partial x_k} \left( v_i \left( p \delta_{ik} + \sigma_{ik}  \right) + \rho T \beta_k \right) = 0. \label{eqn.cl_E}	
\end{equation}
Here, $\x \in \mathbb{R}^d$ is the spatial position vector in $d$ space dimensions, $t \in \mathbb{R}^+$ is the time coordinate and $\A = \left\{ A_{ik} \right\}$ with $A_{ik} = \partial \xi_i / \partial x_k$ is the Eulerian gradient of the Lagrangian position variable $\xi_i$, representing the reference configuration. Furthermore, $\rho$ is the fluid density and $\tau=1/\rho$ is the specific volume, $S$ is the specific entropy, $\vv = \left\{ v_i \right\} \in \mathbb{R}^d$ is the velocity vector, $p$ denotes the fluid pressure, $\bs=\left\{ \sigma_{ik} \right\}$ is the stress tensor, $\J = \left\{ J_i \right\}$ is the specific thermal impulse, $\bG = \left\{ \Gamma_{ik} \right\}$ and $\bbeta = \left\{ \beta_k \right\}$ are thermodynamic dual variables that are defined later and $\ee=E_1 + E_2 + E_3 + E_4$ is the specific total energy, which is a sum of the following contributions: the specific internal energy $E_1 = E_1(\tau,S)$; the kinetic energy $E_2 = \halb v_m v_m$; the energy contained in the elastic deformations $E_3 = \frac{1}{4} c_s^2 \mathring{G}_{ij} \mathring{G}_{ij}$, with the metric tensor $G_{ij} = A_{ki}A_{kj}$, its trace-free part $\mathring{G}_{ij} = G_{ij} - \frac{1}{3} G_{mm} \delta_{ij}$ and $\delta_{ij}$ the usual Kronecker symbol; $E_4 = \halb c_h^2 J_m J_m$ is the energy contained in the specific thermal impulse $J_m$, which is needed within a hyperbolic model of heat conduction. The constant $c_s$ represents the shear sound speed in the medium and $c_h$ is related to the  finite propagation speed of heat waves. The functions $\theta_1>0$ and $\theta_2>0$ depend on two characteristic relaxation times $\tau_1>0$ and $\tau_2>0$ and will be specified later. We recall that Lagrangian mass conservation writes
\begin{equation}
	\frac{\rho}{|\mathbf{A}|} = \rho_0,
	\label{eqn.detAconstraint}
\end{equation}	
where $\rho_0$ is the initial mass density. The equation of the specific volume \eqref{eqn.cl_tau} is redundant since it might be obtained as the dot product of the equation of the distortion tensor \eqref{eqn.cl_A} with the comatrix $|\mathbf{A}|\mathbf{A}^{-1}$ and using the above mass conservation \eqref{eqn.detAconstraint} \cite{BoscheriGPRGCL}.
Furthermore, in the updated Lagrangian framework, the system is supplemented with the trajectory equation for the Eulerian coordinate as 
\begin{equation}
	\mde{\x}{t} = \vv, \qquad \x(0)=\Lx,
	\label{eqn.trajODE}
\end{equation}
 where the Lagrangian coordinate $\boldsymbol{\xi}$ of the initial reference configuration corresponds to the Eulerian coordinate $\x$ at $t=0$. The system is closed by an equation of state (EOS) for the specific internal energy $E_1$ in terms of the specific volume $\tau$ and the specific entropy $\S$, hence $\varepsilon = \varepsilon(\tau,\S)$. In the following we assume that $E_1(\tau,\S)$ is convex in order to assure thermodynamic stability of the system, see \cite{Menikoff198975}. Following \cite{God1961,Rom1998,GPRmodel,HTCGPR,HTCLagrange} the thermodynamic \textit{dual variables}, or \textit{Godunov variables} \cite{Freistuehler2019}, i.e. the \textit{main field} according to \cite{Ruggeri81} read as follows:   
\begin{equation}
	p = - \frac{\partial E}{\partial \tau}, \qquad 
	\vv = \frac{\partial E}{\partial \vv}, \qquad 
	T = \frac{\partial E}{\partial S}, \qquad 
	\bG = \frac{\partial E}{\partial \A}, \qquad 
	\bbeta = \frac{\partial E}{\partial \J}.
\end{equation} 
For convenience we also introduce the vector of dual variables $\w = (-p, \vv, T, \bG, \bbeta)$. 
Based on these definitions, the stress tensor can be expressed as follows: 
\begin{equation}
	\bs = \rho \A^T \bG + \rho \J \otimes \bbeta, \qquad
	\textnormal{or, in index notation} \qquad 
	\sigma_{ik} = \rho A_{ki} \Gamma_{kj} + \rho J_i \beta_j. 
\end{equation}
We furthermore assume that the temperature is positive, $T>0$. After some calculations one can check that the total energy conservation law \eqref{eqn.cl_E} can be obtained as a consequence of \eqref{eqn.PDE_entropy} after multiplication of the governing PDE system with the main field variables $\w$ and subsequent summation.  
In the next section, we design a numerical method that can preserve this thermodynamic compatibility also on the discrete level. 

\section{Numerical method} \label{sec.method}
The computational domain $\Omega(t)$ with boundary $\partial \Omega(t)$ in $d=2$ space dimensions is discretized with a set of non-overlapping triangular control volumes $\omega_c(t)$ of volume $|\omega_c(t)|$. 
\begin{figure}[!htbp]
	\begin{center}
		\begin{tabular}{lll}        
			\includegraphics[trim=50 0 0 50,clip,width=0.40\textwidth]{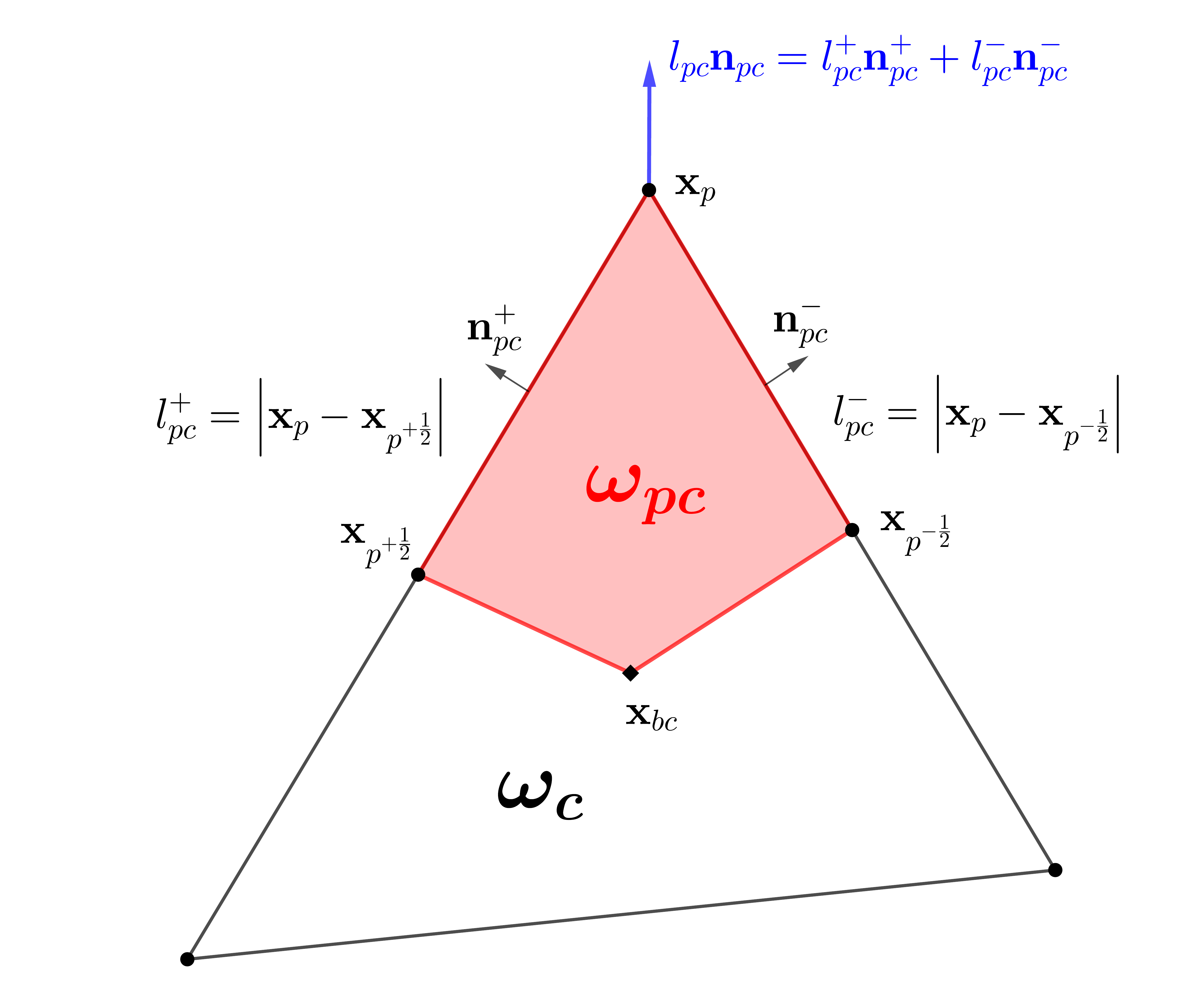} &
			\includegraphics[trim=100 100 100 100,clip,width=0.40\textwidth]{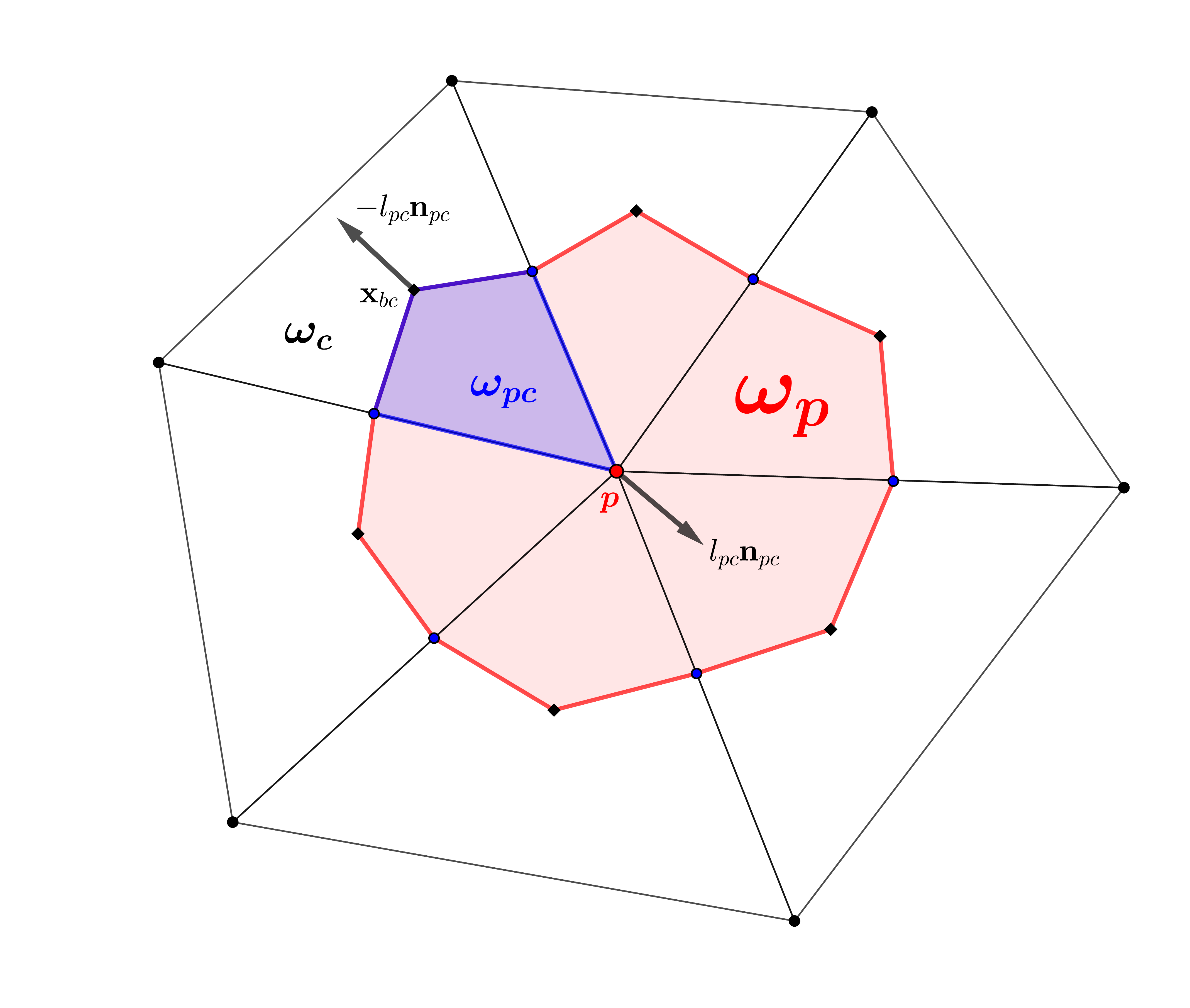} &
		\end{tabular} 
		\caption{Notation for the cell $\omega_c$ (left) and definition of the dual cell $\omega_p$ (right).}
		\label{fig.notation}
	\end{center}
\end{figure}
The quantities related to the mesh geometry are time dependent since the mesh moves, but the topology of the grid remains unchanged. A generic vertex of the mesh is denoted by $p$, while the sub-cell $\omega_{pc}$ refers to the portion of cell $\omega_c$ attached to one of its vertex $p$, as depicted in Figure \ref{fig.notation}. The set of vertexes belonging to cell $\omega_c$ is referred to with $\mathcal{P}(c)$, while the set of cells sharing node $p$ is indicated with $\mathcal{C}(p)$. The sub-cell $\omega_{pc}$ is defined by connecting the vertex with position $\x_p$, the cell barycenter of coordinates $\x_{bc}=1/(d+1) \, \sum_{p \in \mathcal{P}(c)}\x_p$, and the left and right midpoints of the edges impinging on node $p$, namely $p^{-\halb}$ and $p^{+\halb}$, respectively. The half lengths of these edges are given by $l_{pc}^-=\left| \x_p - \x_{p^{-\halb}} \right|$ and $l_{pc}^+=\left| \x_p - \x_{p^{+\halb}} \right|$, and the corner normal is computed as
\begin{equation}
	\lnpc = l_{pc}^+ \n_{pc}^+ + l_{pc}^- \n_{pc}^-.
	\label{eqn.lpcnpc}
\end{equation}
Thanks to Gauss theorem the corner normals satisfy the identity
\begin{equation}
	\sum \limits_{p \in \mathcal{P}(c)} \lnpc = \mathbf{0}.
	\label{eqn.gauss}
\end{equation}
In Figure \ref{fig.notation} we also introduce the dual cell $\omega_p$ given by the union of all sub-cells sharing a generic node $p$. 
The outward pointing corner vectors of $\omega_p$ are simply the corner normals $\lnpc$ with opposite sign, hence also the following identity holds
\begin{equation}
	\sum \limits_{c \in \mathcal{C}(p)} -\lnpc = \mathbf{0}.
	\label{eqn.gaussd}
\end{equation}

\subsection{Semi-discrete thermodynamically compatible Lagrangian scheme} \label{ssec.EC}
In the following we propose a cell-centered finite volume scheme  that is by construction exactly compatible with the entropy inequality \eqref{eqn.cl_S} and which satisfies the extra conservation law of total energy conservation \eqref{eqn.cl_E}. Let $m_c=|\omega_c| \rho_c=|\omega_c| \tau_c^{-1}$ be the mass of the cell, which remains constant in the Lagrangian framework, and let the mass-averaged cell value of a generic quantity $q(\x,t)$ be defined as
\begin{equation}
	q_c = \frac{1}{m_c} \, \int \limits_{\omega_c(t)} \rho \, q \, \text{d}\x.
\end{equation}
Integration of system \eqref{eqn.PDE_entropy} over the cell $\omega_c(t)$, using the Reynolds transport formula and suitable numerical fluctuations in combination with \eqref{eqn.gauss}, yields the following semi-discrete 
Lagrangian HTC scheme:  
\begin{subequations}
	\label{eqn.SD_htc}
	\begin{align}
		&m_c \mde{\tau_c}{t} - \sum \limits_{p \in \mathcal{P}(c)} \lnpc \cdot ({\vv}_p-\vv_c) = 0,    \label{eqn.sdhtc_tau} \\
		&m_c \mde{\vv_c}{t} + \sum \limits_{p \in \mathcal{P}(c)} \left( \lnpc \, ({p}_p-p_c)  +  \lnpc \cdot (\bs_p-\bs_c) \right) \label{eqn.sdhtc_v} \\
		&  \phantom{m_c \mde{\vv_c}{t}} + \sum \limits_{p \in \mathcal{P}(c)} l_{pc} (\alpha_p + \epsilon_p ) \, (\vv_c- \vv_p)  = \mathbf{0},    \nonumber  \\ 
		&m_c \mde{\S_c}{t} +  \sum \limits_{p \in \mathcal{P}(c)} \lnpc \cdot  (\rho_p \bbeta_p- \rho_c \bbeta_c) = \Pi_c + \pi_c \geq 0.  \label{eqn.sdhtc_S} \\ 
		&m_c \mde{\A_c}{t}   + \rho_c \A_c \sum \limits_{p \in \mathcal{P}(c)} \lnpc   (\vv_p-\vv_c) = -\frac{1}{\theta_1(\tau_1)} \boldsymbol{\Gamma}_c, \label{eqn.sdhtc_A}  \\ 		
		&m_c \mde{\J_c}{t}  		
		+ \rho_c \J_c \cdot \sum \limits_{p \in \mathcal{P}(c)} \lnpc   (\vv_p-\vv_c) 
		+ \rho_c \sum \limits_{p \in \mathcal{P}(c)} \lnpc \,  (T_p-T_c)
		\label{eqn.sdhtc_J} 
		 \\
		& \phantom{m_c \mde{\J_c}{t}} 
		+ \rho_c \sum \limits_{p \in \mathcal{P}(c)} \left( \lnpc  \, (\alpha_p + \epsilon_p ) |\omega_p| \nabla_p^c \cdot (\rho_c \bbeta_c)  \right) = -\frac{1}{\theta_2(\tau_2)} \bbeta_c,  \nonumber   
	\end{align}
\end{subequations} 
with $\sum_p \lnpc \cdot \boldsymbol{\sigma}_p = 
\sum_p l_{pc} n_{pc}^{k} \cdot {\sigma}^{ik}_p$ and $\A_c \sum_p  \lnpc   \tilde{\vv}_p = 
A^{il}_c \sum_p l_{pc} n_{pc}^{k}   \tilde{v}^{l}_p$,
where $n_{pc}^{k}$ is the $k$-th component of the unit normal vector $\mathbf{n}_{pc}$.
As in \cite{HTCLagrange} we have followed the general framework of entropy conservative schemes proposed by Abgrall \textit{et al.} in \cite{Abgrall2018,HTCAbgrall}, with the scalar correction factor $\alpha_p \in \mathbb{R}$ that ensures discrete thermodynamic compatibility with the total energy conservation law and which will be defined later. $\epsilon_p \geq 0$ is a numerical viscosity that can be switched off to obtain an entropy conserving Lagrangian scheme (ECL). The object $|\omega_p| \nabla_p^c \cdot \rho_c \bbeta_c$ is a volume-integrated discrete divergence defined for each node $p$ as 
\begin{equation}
	|\omega_p| \nabla_p^c \cdot \rho_c \bbeta_c = -\sum \limits_{c \in \mathcal{C}(p)} \lnpc \cdot \rho_c \bbeta_c.
	\label{eqn.divrhobeta}
\end{equation}
Furthermore, $\rho_p$, $p_p$, $T_p$ are the averaged density, pressure and temperature in the node $p$, while   ${\boldsymbol{\sigma}}_p$ and  $\bbeta_p$ are node-averages of the stress tensor and of the thermodynamic dual variable to the thermal impulse, respectively. 
The mesh motion is governed by the trajectory equation \eqref{eqn.trajODE} that is semi-discretized at each node of the mesh as 
\begin{equation}
	\mde{\x_p}{t} = \vv_p, \qquad \x_p(0)= \Lx_p(0).
	\label{eqn.xp}
\end{equation}
The nodal values of $\rho_p$, $\vv_p$, $p_p$, $T_p$, $\bbeta_p$ and $\bs_p$ are obtained via conservation principles, i.e. requiring that the sum of the fluctuations around a node is equal to the sum of the fluxes through the boundary of the dual cell $\omega_p$, hence, to achieve conservation in the discrete momentum equation \eqref{eqn.sdhtc_v} we require
\begin{eqnarray}
	\label{eqn.conservation}
	&& \phantom{-} \sum \limits_{c \in \mathcal{C}(p)} \lnpc \, (p_p-p_c) + 
	\sum \limits_{p \in \mathcal{P}(c)} \lnpc \cdot (\bs_p-\bs_c)  
	+
	\alpha_p  \sum \limits_{c \in \mathcal{C}(p)} l_{pc} (\vv_c-\vv_p) 
	\\
	&=& - \sum \limits_{c \in \mathcal{C}(p)} \lnpc \, p_c \,  
	    - \sum \limits_{c \in \mathcal{C}(p)} \lnpc \cdot \bs_c \, 
	=  \int \limits_{\partial \omega_p(t)} \left( p \mathbf{I} + \bs\right) \cdot \n \, \text{ds}
	\nonumber .
\end{eqnarray}
Using \eqref{eqn.gaussd}, we obtain the nodal velocity that satisfies conservation: 
\begin{equation}
	\vv_p = \frac{\sum \limits_{c \in \mathcal{C}(p)}  l_{pc} \, \vv_c}{\sum \limits_{c \in \mathcal{C}(p)} l_{pc}}. \qquad 
	\label{eqn.vp}
\end{equation}
For $\rho_p$, $p_p$, $T_p$, $\bbeta_p$ and $\bs_p$ the choice of the corner values is arbitrary, since conservation is automatically guaranteed, so we simply choose the same average as above, i.e.  
\begin{equation}
	p_p = \frac{\sum \limits_{c \in \mathcal{C}(p)} l_{pc} \, p_c}{\sum \limits_{c \in \mathcal{C}(p)} l_{pc}}, \quad 
	\rho_p = \frac{\sum \limits_{c \in \mathcal{C}(p)} l_{pc} \, \rho_c}{\sum \limits_{c \in \mathcal{C}(p)} l_{pc}}, \quad 
	T_p = \frac{\sum \limits_{c \in \mathcal{C}(p)} l_{pc} \, T_c}{\sum \limits_{c \in \mathcal{C}(p)} l_{pc}},
	\label{eqn.prhosT.1}
\end{equation}
\begin{equation}
	\bs_p = \frac{\sum \limits_{c \in \mathcal{C}(p)}  l_{pc} \, \bs_c}{\sum \limits_{c \in \mathcal{C}(p)} l_{pc}}, \quad 
	\bbeta_p = \frac{\sum \limits_{c \in \mathcal{C}(p)}  l_{pc} \, \bbeta_c}{\sum \limits_{c \in \mathcal{C}(p)} l_{pc}}.		
	\label{eqn.prhosT.2}
\end{equation}

As done in \cite{HTCLagrange}, we now require thermodynamic compatibility by imposing the following condition around a node $p$: the sum of the total energy fluctuations around a node must be equal to the sum of total energy fluxes around the same node. The total energy fluctuations are obtained as the dot product of the fluctuations in the scheme \eqref{eqn.SD_htc} with the main-field variables $\w_c$. This leads to the following nodal condition: 
\begin{eqnarray}
	\qquad &&  \phantom{+} \sum \limits_{c \in \mathcal{C}(p)} \! \lnpc \cdot \left( p_c \, (\vv_p-\vv_c) + \vv_c \, (p_p-p_c) +  (\bs_p-\bs_c) \, \vv_c + T_c \, (\rho_p \bbeta_p - \rho_c \bbeta_c) \right)  	\label{eqn.Ehtc}
	 \\ 
	&&  + \sum \limits_{c \in \mathcal{C}(p)} \! \lnpc \cdot \left( 
	\rho_c \bG_c \A_c (\vv_p - \vv_c) +  \rho_c \bbeta_c \otimes \J_c (\vv_p - \vv_c)
	+ 	\rho_c \bbeta_c \, (T_p-T_c)  \right) 	
	\nonumber \\ 
&&  + \alpha_p \left( \sum \limits_{c \in \mathcal{C}(p)} l_{pc} \, \vv_c \cdot (\vv_c-\vv_p) + \sum \limits_{c \in \mathcal{C}(p)} \lnpc \cdot \rho_c \bbeta_c ( |\omega_p| \nabla_p^c \cdot \rho_c \bbeta_c ) \right) \nonumber \\ 
	&&  = - \sum \limits_{c \in \mathcal{C}(p)} \! \lnpc \cdot \left(  p_c \, \vv_c + \bs_c \vv_c + \rho_c T_c \bbeta_c \right) = \int \limits_{\partial \omega_p(t)} \left( p \, \vv + \bs \vv + \rho T \bbeta \right) \cdot \n \, \text{ds}. \nonumber 
\end{eqnarray}
The terms multiplied by the scalar factor $\alpha_p$ can be conveniently rearranged in quadratic form using the definitions \eqref{eqn.vp} and \eqref{eqn.divrhobeta}, see also
\cite{HTCLagrange}:
\begin{eqnarray}
	\label{eqn.alphaterms2}	
	  \delta_p &:= & \sum \limits_{c \in \mathcal{C}(p)} l_{pc} \, \vv_c \cdot (\vv_c-\vv_p) + \sum \limits_{c \in \mathcal{C}(p)} \lnpc \cdot \rho_c \bbeta_c ( |\omega_p| \nabla_p^c \cdot \rho_c \bbeta_c )  \nonumber \\ 
	  & = & 
	   \left( \sum \limits_{c \in \mathcal{C}(p)} l_{pc} \, (\vv_c-\vv_p)^2 \right)  +
	  \left( |\omega_p| \nabla_p^c \cdot \rho_c \bbeta_c \right)^2 \geq 0.   
\end{eqnarray}
Solving \eqref{eqn.Ehtc} for $\alpha_p$ one obtains  
\begin{equation}
\alpha_p = \frac{\nu_p}{\delta_p}, 
\label{eqn.alpha} 
\end{equation}
with the denominator already given in \eqref{eqn.alphaterms2} and the numerator  
\begin{eqnarray}
	\nu_p &=& \phantom{+} \sum \limits_{c \in \mathcal{C}(p)} \!\!\! \lnpc \cdot (p_c \, \vv_c - p_c \, \vv_p  - p_p \, \vv_c) 
	+ 
	\sum \limits_{c \in \mathcal{C}(p)} \!\!\! \lnpc \cdot (\bs_c \vv_c - \bs_c \vv_p - \bs_p \vv_c) 
	\nonumber \\ 
	&& + 
	\sum \limits_{c \in \mathcal{C}(p)} \!\!\! \lnpc \cdot (T_c \rho_c  \bbeta_c - T_c \rho_p \bbeta_p   - T_p  \rho_c \bbeta_c).
	\label{eqn.nu}
\end{eqnarray}
For $\delta_p=0$ we set $\alpha_p=0$. The two entropy production terms read as follows: the one related to the algebraic source term is simply given by 
\begin{equation}
	\pi_c = \frac{\boldsymbol{\Gamma}_c:\boldsymbol{\Gamma}_c}{T_c \theta_1(\tau_1)}  +
	\frac{\bbeta_c \cdot \bbeta_c}{T_c \theta_2(\tau_2)}  \geq 0, 
\end{equation} 
which is quadratic in $\bG_c$ and $\bbeta_c$ and clearly non-negative since $T_c>0$ and $(\theta_1, \theta_2) > 0$ is assumed. The entropy production term related to the numerical viscosity $\epsilon_p \geq 0$ that may be introduced into the numerical scheme reads  
\begin{equation}
	\Pi_c = \sum \limits_{p \in \mathcal{P}(c)} \frac{|\omega_{pc}|}{|\omega_c|} \frac{\Pi_p}{T_c}, 
\end{equation}
with the nodal numerical entropy production  
\begin{equation} 	
	\Pi_p = \epsilon_p \left( \sum \limits_{c \in \mathcal{C}(p)} l_{pc} \, (\vv_c-\vv_p)^2 \right)  +
	\epsilon_p \left( |\omega_p| \nabla_p^c \cdot \rho_c \bbeta_c \right)^2 = \epsilon_p \delta_p \geq 0, 
\end{equation}
which is also obviously non-negative. Since by assumption $T_c >0$ also $\Pi_c \geq 0$. Hence, in our scheme, we immediately have a discrete cell entropy inequality, \textit{by construction}. 

As we prove now, the semi-discrete scheme \eqref{eqn.SD_htc} is also nonlinearly stable in the sense that it conserves total energy exactly, it satisfies the determinant constraint $|\A_c|=\rho_c / \rho_c(0)$ at a semi-discrete level and preserves the curl-free property of $\J_c$ and $\A_c$ for initially curl-free data and vanishing sources ($\tau_1 \to \infty$ and $\tau_2 \to \infty$). To show the involution-preserving properties we first prove the following Lemma, see also \cite{Barsukow2024,Sidilkover2025,CompatibleDG1}.

\begin{lemma} \label{lemma.rot.grad} 
	In two space dimensions ($d=2$) the discrete gradient operator applied to a scalar field $\phi_p$ defined in the nodes $p$ 
	\begin{equation}
		\nabla_c^p \phi_p = \frac{1}{|\omega_c|} \sum \limits_{p \in \mathcal{P}(c)} \lnpc \phi_p
		\label{eqn.nablacp} 
	\end{equation} 
	and its dual discrete curl operator applied to a vector field $\J_c$ defined in the cell centers 
	\begin{equation}
		\nabla_p^c \times \J_c = - \frac{1}{|\omega_p|} \sum \limits_{c \in \mathcal{C}(p)} \lnpc \times \J_c	
		\label{eqn.curlpc} 
	\end{equation} 
	satisfy the discrete vector calculus identity
	\begin{equation}
		\nabla_p^c \times \nabla_c^p \phi_p = 0. 
		\label{eqn.disc.rot.grad} 
	\end{equation}	
\end{lemma}

\begin{proof}
	In the following we denote the indices of the three nodes of control volume $\omega_c$ by $p_1$, $p_2$ and $p_3$, with Eulerian node coordinates $\x_{p_1}$, $\x_{p_2}$ and $\x_{p_3}$, respectively. We furthermore assume that the nodes are ordered counter-clockwise. The indices of the edges of $\omega_c$ are denoted by $e_1$, $e_2$ and $e_3$, respectively, with edge $e_1$ composed of nodes $p_1$ and $p_2$, edge $e_2$ composed of nodes $p_2$ and $p_3$, and edge $e_3$ composed of nodes $p_3$ and $p_1$, respectively. Assuming a linear distribution of the scalar field $\phi_p$ inside control volume $\omega_c$ (according to continuous P1 Lagrange finite elements) immediately yields the following averages in the edge midpoints: 
	\begin{equation}
		\phi_{e_1} = \halb \left( \phi_{p_1} + \phi_{p_2} \right), \quad 
		\phi_{e_2} = \halb \left( \phi_{p_2} + \phi_{p_3} \right), \quad 
		\phi_{e_3} = \halb \left( \phi_{p_3} + \phi_{p_1} \right). 
	\end{equation}
	It can easily be checked that the three corner normals are given by 
	\begin{eqnarray}
		l_{{p_1}c} \n_{{p_1}c} &=& \halb \left( (y_{p_2} - y_{p_3}),  (x_{p_3} - x_{p_2}), 0 \right)^T, \nonumber \\ 
		l_{{p_2}c} \n_{{p_2}c} &=& \halb \left( (y_{p_3} - y_{p_1}),  (x_{p_1} - x_{p_3}), 0 \right)^T, \nonumber \\ 
		l_{{p_3}c} \n_{{p_3}c} &=& \halb \left( (y_{p_1} - y_{p_2}),  (x_{p_2} - x_{p_1}), 0 \right)^T, 
	\end{eqnarray}
	while the area of element $\omega_c$ is 
	\begin{equation}
		|\omega_c| = \halb \left( x_{p_1} y_{p_2} - x_{p_1} y_{p_3} - x_{p_2} y_{p_1} + x_{p_2} y_{p_3} + x_{p_3} y_{p_1} - x_{p_3} y_{p_2} \right). 
	\end{equation}
	The discrete gradient $\nabla_c^p \phi_p = \frac{1}{|\omega_c|} \sum_p \lnpc \phi_p$ of the scalar field $\phi_p$ in the cell $c$ then reads 
	\begin{equation}
		\nabla_c^p \phi_p = \frac{1}{2 |\omega_c|}
		\left( \begin{array}{c} (y_{p_2} - y_{p_3}) \phi_{p_1} + (y_{p_3} - y_{p_1}) \phi_{p_2} + (y_{p_1} - y_{p_2}) \phi_{p_3} \\
			                    (x_{p_3} - x_{p_2}) \phi_{p_1} + (x_{p_1} - x_{p_3}) \phi_{p_2} + (x_{p_2} - x_{p_1}) \phi_{p_3} \\
			 0
			 \end{array} 
			  \right).
	\end{equation}
	Since this property will be needed later, we now also show that this gradient is identical to the discrete gradient obtained by using continuous P1 Lagrange finite elements inside the control volume $\omega_c$. The associated basis functions in the universal reference element $T_0 = {\left\{ \zeta_1, \zeta_2: 0 \leq \zeta_1 \leq 1, 0 \leq \zeta_2 \leq 1-\zeta_1 \right\}}$ with nodes $(0,0)$, $(1,0)$ and $(0,1)$ and reference coordinates $\boldsymbol{\zeta}= (\zeta_1, \zeta_2)$ read 
	\begin{equation}
		\psi_{p_1} = \psi_{p_1}(\boldsymbol{\zeta}) = 1- \zeta_1 - \zeta_2, \qquad 
		\psi_{p_2} = \psi_{p_2}(\boldsymbol{\zeta}) = \zeta_1, \qquad 
		\psi_{p_3} = \psi_{p_3}(\boldsymbol{\zeta}) = \zeta_2. 
	\end{equation}
	The Jacobian of the mapping 
	\begin{equation}
		\x = \x(\boldsymbol{\zeta}) = 
		\psi_{p_1}(\boldsymbol{\zeta}) \, \x_{p_1} + 
		\psi_{p_2}(\boldsymbol{\zeta}) \, \x_{p_2} + 
		\psi_{p_3}(\boldsymbol{\zeta}) \, \x_{p_3} 
		= \sum \limits_{p \in \mathcal{P}(c)} \psi_p \x_p . 
	\end{equation}
	reads
	\begin{equation}
		\mathbf{M}_c(\x) = \frac{\partial \x}{\partial \boldsymbol{\zeta}} = 
		\left( \begin{array}{cc}
				   x_{p_2} - x_{p_1} & x_{p_3} - x_{p_1} \\ 
				   y_{p_2} - y_{p_1} & y_{p_3} - y_{p_1}
			   \end{array}
		\right) = 
		\sum \limits_{p \in \mathcal{P}(c)}  \frac{\partial \psi_p}{\partial \boldsymbol{\zeta}} \, \x_p  ,
		\label{eqn.Mc}
	\end{equation}
	with $|\mathbf{M}_c| = 2 |\omega_c|$. 
	The discrete gradient of a discrete scalar quantity $\phi_h = \phi_{p_1} \psi_{p_1} + \phi_{p_2} \psi_{p_2} + \phi_{p_3} \psi_{p_3}$ is then given by 
	\begin{eqnarray}
		 \tilde{\nabla}_c^p \phi_p &=& \mathbf{M}_c^{-T} \frac{\partial \phi_h}{\partial \boldsymbol{\zeta} }
		  = \mathbf{M}_c^{-T} \sum \limits_{p \in \mathcal{P}(c)}  \, \frac{\partial \psi_p}{\partial \boldsymbol{\zeta} }
		  \, 
		  \phi_p		 
		  \nonumber \\ 
		 &=& \frac{1}{2 |\omega_c|}
		 \left( \begin{array}{c} (y_{p_2} - y_{p_3}) \phi_{p_1} + (y_{p_3} - y_{p_1}) \phi_{p_2} + (y_{p_1} - y_{p_2}) \phi_{p_3} \\
		 	(x_{p_3} - x_{p_2}) \phi_{p_1} + (x_{p_1} - x_{p_3}) \phi_{p_2} + (x_{p_2} - x_{p_1}) \phi_{p_3} \\
		 	0
		 \end{array} 
		 \right) = \nabla_c^p \phi_p. 
	\end{eqnarray}	
	In the following we will denote a generic node of element $\omega_c$ by $p$, as usual, and the corresponding two nodes on the opposite edge are denoted in counter-clockwise order by $p^-$ and $p^+$, respectively. Furthermore, we will denote the edge that connects $p$ with $p^-$ by $e_{pc}^-$ and the edge that links $p$ with $p^+$ by $e_{pc}^+$. 
	We now compute the single contribution of element $\omega_c$ and node $p$ to the curl on the dual mesh $|\omega_{p}| \nabla_{p}^c \times \nabla_c^{p} \phi_p = - \sum_c \lnpc \times \nabla_c^p \phi_p$, i.e. the contribution $-\lnpc \times \nabla_c^p \phi_p$, which after a lengthy but straightforward calculation reduces to 
	\begin{equation}
		-\lnpc \times \nabla_c^p \phi_p = \halb \left( \phi_{p^+} - \phi_{p^-} \right) = \phi_{e_{pc}^+} - \phi_{e_{pc}^-},
	\end{equation}  
	i.e. the difference of values at the respective edge midpoints. 
	Consider now two adjacent elements $\omega_c$ and $\omega_d$ which share a common edge $e$ composed of nodes $p$ and $q$. We assume that the elements are ordered counter clockwise, i.e. the previous element is $\omega_c$ and the subsequent element is $\omega_d$. The value at the edge midpoint given by the previous element $\omega_c$ is denoted by $\phi_e^- = \phi_{e_{pc}^+}$ and the value at the edge midpoint given by the subsequent element $\omega_d$ is denoted by $\phi_e^+ = \phi_{e_{pd}^-}$, 
	Since the values of the scalar field $\phi$ at the edge midpoint are the \textit{same} for two adjacent elements  that share a common edge $e$, i.e. $\phi_e^+ = \phi_e^- = \halb (\phi_p + \phi_q) $, summing up over all elements $\omega_c$ around the node $p$ yields zero, due to the telescopic sum property:
 	\begin{equation}
 		|\omega_{p}| \nabla_{p}^c \times \nabla_c^{p} \phi_p = - \sum \limits_{c \in \mathcal{C}(p)} \lnpc \times \nabla_c^p \phi_p = \sum \limits_{e \in \mathcal{E}(p)} \phi_e^+ - \phi_e^- = 0. 	
 	\end{equation}
 	Here, $\mathcal{E}(p)$ is the set of edges connected to node $p$. The above identity completes the proof that the discrete curl on the dual mesh applied to the discrete gradient on the primal mesh is zero. 
\end{proof}
This Lemma has already been proven previously in a different setting in \cite{Sidilkover2025} and \cite{CompatibleDG1}.

\begin{theorem} \label{th.NLstability}
Assuming periodic boundary conditions or the boundary conditions $\int \limits_{\partial \Omega} \vv \cdot \n \, \textnormal{ds}= 0$, $\int \limits_{\partial \Omega} \bs \cdot \n \, \textnormal{ds}= 0$ and $\int \limits_{\partial \Omega} \bbeta \cdot \n \, \textnormal{ds}= 0$,
the semi-discrete scheme \eqref{eqn.SD_htc} with node values given by \eqref{eqn.vp}, \eqref{eqn.prhosT.1} and \eqref{eqn.prhosT.2} and the correction factor given by
\eqref{eqn.Ehtc} and \eqref{eqn.alpha} is nonlinearly stable in the sense that it conserves total energy: 
\begin{equation}
	\int \limits_{\Omega} \mde{\ee}{t} \, \text{d}\x = 0.
	\label{eqn.NLstability}
\end{equation}	
\end{theorem}

\begin{proof}
The semi-discrete total energy conservation law can be obtained by contracting the semi-discrete scheme \eqref{eqn.SD_htc} with the discrete main field  $\w_c=(-p_c,\vv_c,T_c,\bG_c,\bbeta_c)^T$. 
Summing up over all $c$ and exchanging summation over $p$ and $c$ in the fluctuation terms we thus obtain 
\begin{eqnarray}
	&  \sum_c m_c \left( -p_c  \mde{\tau_c}{t} + \vv_c \mde{\vv_c}{t} 
	   + T_c \mde{\S_c}{t} + \bG_c  : \mde{\A_c}{t} + 
	   \bbeta_c  \cdot \mde{\J_c}{t} \right)  	    
	   \nonumber \\ 	
	& 
	 + \sum_p \sum \limits_{c \in \mathcal{C}(p)} \lnpc \cdot \left( p_c ({\vv}_p-\vv_c)  +
	  \vv_c ({p}_p-p_c)  +   (\bs_p-\bs_c) \vv_c + T_c (\rho_p \bbeta_p- \rho_c \bbeta_c)  \right) \nonumber 
	    \nonumber \\ 
	&   + \sum_p \sum \limits_{c \in \mathcal{C}(p)} \lnpc \cdot \left( \rho_c \bG_c \A_c   (\vv_p-\vv_c) + \rho_c   \bbeta_c \otimes \J_c  (\vv_p-\vv_c)  + \rho_c \bbeta_c  (T_p-T_c)  \right) \nonumber \\ 
	& 
	+ \sum_p \sum \limits_{c \in \mathcal{C}(p)} l_{pc} (\alpha_p + \epsilon_p ) \, (\vv_c- \vv_p) \cdot \vv_c 
     \nonumber \\ 
	& 
	+  \sum_p \sum \limits_{c \in \mathcal{C}(p)}  \lnpc  \cdot \rho_c \bbeta_c \, (\alpha_p + \epsilon_p ) |\omega_p| \nabla_p^c \cdot (\rho_c \bbeta_c)   = 
	\nonumber \\
	& 
	\sum_c \left( T_c \Pi_c  + T_c \pi_c - \frac{1}{\theta_2(\tau_2)} \boldsymbol{\beta}_c \cdot \boldsymbol{\beta}_c - \frac{1}{\theta_1(\tau_1)} \boldsymbol{\Gamma}_c : \boldsymbol{\Gamma}_c \right). 
  \nonumber   
\end{eqnarray} 
The term in the first bracket is obviously the time derivative of the total energy in cell $c$. Thanks to the definition of $\alpha_p$ given in \eqref{eqn.alpha} and \eqref{eqn.Ehtc}, using the boundary conditions and realizing that the sum over all $p$ of the boundary integrals on the dual volumes shown in \eqref{eqn.Ehtc} cancel out 
and making use of the definition of the entropy production terms $\Pi_c$ and $\pi_c$, the above expression reduces to 
\begin{equation}
	\sum_c m_c  \mde{E_c}{t}  = \int \limits_{\Omega} \mde{\ee}{t} \, \text{d}\x = 0, 
\end{equation} 
which completes the proof.  	    
 
\end{proof}		

\begin{theorem} \label{th.detA}
	The semi-discrete scheme \eqref{eqn.SD_htc} preserves the determinant constraint $|\A_c|=\rho_c(t)/\rho_c(0)$ exactly on the semi-discrete level in the sense  
	\begin{equation}
	\frac{d|\A_c|}{dt} = \frac{d}{dt} \left( \frac{\rho_c}{\rho_c(0)} \right) = - \frac{\rho_c^2}{\rho_c(0)} \frac{d \tau_c}{dt} = -\frac{|\A_c|}{|\omega_c|} \sum \limits_{p \in \mathcal{P}(c)} \lnpc \cdot ({\vv}_p-\vv_c).
		\label{eqn.detA.pres}
	\end{equation}	
\end{theorem}

\begin{proof}		
	Division of \eqref{eqn.sdhtc_A} by $m_c = \rho_c |\omega_c| $ leads to  
	\begin{equation}	
	\mde{\A_c}{t}   = -\frac{\A_c}{|\omega_c|} \sum \limits_{p \in \mathcal{P}(c)} \lnpc   (\vv_p-\vv_c).
	\end{equation} 
	Using the Jacobi identity for the time derivative of the determinant of $\A_c$ and using the discretized PDE for $\tau_c = \rho_c^{-1}$ given by \eqref{eqn.sdhtc_tau} yields the sought result, which is the discrete counterpart of \eqref{eqn.detAconstraint}:   
	\begin{eqnarray}	
	\mde{|\A_c|}{t}   &=& \phantom{-} |\A_c| \tr \left( \A^{-1} \mde{\A_c}{t} \right) \nonumber \\ 
	&=& -|\A_c| \tr \left( \A_c^{-1} \A_c 
	\frac{1}{|\omega_c|} \sum \limits_{p \in \mathcal{P}(c)} \lnpc   (\vv_p-\vv_c) \right) 
	\nonumber \\ 
	&=& -\frac{ |\A_c|}{|\omega_c|} \sum \limits_{p \in \mathcal{P}(c)} \lnpc \cdot  (\vv_p-\vv_c) = - \frac{\rho_c^2}{\rho_c(0)} \frac{d \tau_c}{dt} = \frac{d}{dt} \left( \frac{\rho_c}{\rho_c(0)} \right).
\end{eqnarray} 
	
\end{proof}

\begin{theorem} \label{th.curlA}
	For vanishing source, i.e. $\tau_1 \to \infty$, and initially curl-free data for the distortion field given by $\A_c(0) = \nabla_c^p \Lx_p$, hence $\nabla_p^c \times \A_c(0) = 0$, and assuming exact time integration, the semi-discrete scheme \eqref{eqn.SD_htc} preserves the curl of the rows of $\A_c$ exactly at the semi-discrete level in the sense 
	\begin{equation}
		\nabla_p^c \times \A_c(t) = 0.
		\label{eqn.curlA.pres}
	\end{equation}	
\end{theorem}

\begin{proof}		
	In the absence of source terms ($\tau_1 \to \infty$) the distortion field represents the inverse deformation gradient, i.e. the Eulerian gradient of the Lagrangian reference coordinates: 
	\begin{equation}
		 \A_c = \nabla_c^p \Lx_p = \left( 
		 \sum \limits_{p \in \mathcal{P}(c)}  \, \frac{\partial \psi_p}{\partial \boldsymbol{\zeta} }
		 \, 
		 \Lx_p \right) 
		 \mathbf{M}_c^{-1},
		 \label{eqn.A.def}   
	\end{equation} 
	which is in particular true for the initial time $t=0$. 
	Division of \eqref{eqn.sdhtc_A} by $m_c = \rho_c |\omega_c|$ and using  $\sum_p \lnpc \vv_c = 0$ yields  
	\begin{equation} 
	\mde{\A_c}{t}  		
	+  \A_c  \, \frac{1}{\omega_c} \sum \limits_{p \in \mathcal{P}(c)} \lnpc  \vv_p  = 0, 
	\end{equation} 	
	or, in more compact notation 
	\begin{equation} 
		\mde{\A_c}{t}  		
		+  \A_c \, \nabla_c^p  \vv_p = 0. 
		\label{eqn.A.simple}
	\end{equation} 
	Applying the Lagrangian time derivative to the definition \eqref{eqn.A.def} yields 
	\begin{equation}
	 \mde{\nabla_c^p \Lx_p}{t} =
	 \left( 
	  \sum \limits_{p \in \mathcal{P}(c)}  \, \frac{\partial \psi_p}{\partial \boldsymbol{\zeta} }
	\, 	\Lx_p \right)
	\mde{\mathbf{M}_c^{-1}}{t} 
	\nonumber \\
	= - \left( \sum \limits_{p \in \mathcal{P}(c)}  \, \frac{\partial \psi_p}{\partial \boldsymbol{\zeta} }
	\, 	\Lx_p \right) 
	\mathbf{M}_c^{-1} \, \frac{d \mathbf{M}_c}{dt} \, \mathbf{M}_c^{-1}
	,	
	\label{eqn.dAdtaux} 
	\end{equation} 	
	where we have employed the formula for the time derivative of the inverse of a matrix and where we have also used the fact that the only object which depends on time is actually $\mathbf{M}_c$ due to the moving mesh. Indeed, at the aid of the discrete trajectory equation \eqref{eqn.xp} and thanks to \eqref{eqn.Mc} we obtain the following time derivative of $\mathbf{M}_c$: 
	\begin{equation} 
	\frac{d \mathbf{M}_c}{dt} = \sum \limits_{p \in \mathcal{P}(c)} \frac{\partial \psi_p}{\partial \boldsymbol{\zeta}} \, \frac{d\x_p}{dt} = 
	\sum \limits_{p \in \mathcal{P}(c)} \frac{\partial \psi_p}{\partial \boldsymbol{\zeta}} \, \vv_p. 
	\label{eqn.dMcdt} 
	\end{equation}
	Inserting \eqref{eqn.dMcdt} into \eqref{eqn.dAdtaux} and using \eqref{eqn.A.def} leads to 
	\begin{equation}
	\mde{\nabla_c^p \Lx_p}{t} = - \left( \sum \limits_{p \in \mathcal{P}(c)} \frac{\partial \psi_p}{\partial \boldsymbol{\zeta}} \, \Lx_p  \mathbf{M}_c^{-1} \right) \, \left(  \sum \limits_{q \in \mathcal{P}(c)}  \, \frac{\partial \psi_q}{\partial \boldsymbol{\zeta} }
	\, 	\vv_q \, \mathbf{M}_c^{-1} \right) = - \nabla_c^p \Lx_p \, \nabla_c^q  \vv_q,	
	\end{equation} 
	which corresponds to \eqref{eqn.A.simple} and thus to the numerical scheme \eqref{eqn.sdhtc_A} used for the discretization of $\A_c$.  We thus conclude that \eqref{eqn.sdhtc_A} is equivalent to 
	\begin{equation}
		 \mde{\A_c}{t} = \mde{\nabla_c^p \Lx_p}{t}. 
	\end{equation}
	This means that, assuming exact time integration, $\A_c = \nabla_c^p \Lx_p$ for all times, since this relation is assumed to be true at the initial time $\A_c(0) = \nabla_c^p \Lx_p$. Application of the discrete curl operator on the dual mesh to the previous relation leads to 
	\begin{equation}
	 \nabla_p^c \times \A_c  = 
	 \nabla_p^c \times \nabla_c^q \Lx_q  = 0, 
	\end{equation}
	thanks to the discrete vector calculus identity \eqref{eqn.disc.rot.grad}, which completes the proof. 	 
\end{proof}

\begin{theorem} \label{th.curlJ}
	For vanishing source, i.e. $\tau_2 \to \infty$, and initially curl-free data for the thermal impulse given by the discrete gradient of a scalar potential $Z$ as 
	$\J_c(0) = \nabla_c^p Z_p(0)$ and which thus satisfies $\nabla_p^c \times \J_c(0) = 0$, assuming exact time integration, the semi-discrete scheme \eqref{eqn.SD_htc} preserves the curl-free condition of $\J_c$ exactly at the semi-discrete level for all times, i.e. 
	\begin{equation}
	\nabla_p^c \times \J_c = 0.
		\label{eqn.curlJ.pres}
	\end{equation}	
\end{theorem}

\begin{proof}			
	Division of \eqref{eqn.sdhtc_J} by $m_c = \rho_c |\omega_c|$, introducing the auxiliary scalar field $\auxf_p = T_p + \alpha_p |\omega_p| \nabla_p^c \cdot (\rho_c \bbeta_c) $ and using the fact that $\sum_p \lnpc T_c = 0$ and $\sum_p \lnpc \vv_c = 0$ yields 
	\begin{equation} 
	\mde{\J_c}{t}  		
	+  \J_c \cdot \frac{1}{\omega_c} \sum \limits_{p \in \mathcal{P}(c)} \lnpc  \vv_p
	+  \frac{1}{\omega_c} \sum \limits_{p \in \mathcal{P}(c)} \lnpc \auxf_p  = 0, 
	\end{equation} 
	or, in more compact notation 
	\begin{equation} 
	\mde{\J_c}{t}  		
	+  \J_c \cdot \nabla_c^p  \vv_p
	+  \nabla_c^p \auxf_p = 0. 
	\label{eqn.J.compact}
	\end{equation} 
 	We now define the discrete gradient of a discrete scalar potential $Z$ 
 		\begin{equation}
 		\J_c = \nabla_c^p Z_p = \left( 
 		\sum \limits_{p \in \mathcal{P}(c)}  \, \frac{\partial \psi_p}{\partial \boldsymbol{\zeta} }
 		\, 
 		Z_p(t) \right) \cdot 
 		\mathbf{M}_c^{-1},
 		\label{eqn.J.def}   
 	\end{equation}  
 	which is in particular assumed to hold at the initial time $t=0$.  
 	We furthermore define 
 	\begin{equation} 
 		\frac{d Z_p}{dt} = -\auxf_p.
 	\end{equation} 
	In the relation \eqref{eqn.J.def}, both $Z_p$ and $\M_c$ are functions of time. 
	Applying the Lagrangian time derivative yields 
	\begin{eqnarray}
		\frac{d \nabla_c^p Z_p}{dt} &=& \left( 
		\sum \limits_{p \in \mathcal{P}(c)}  \, \frac{\partial \psi_p}{\partial \boldsymbol{\zeta} }
		\, 
		Z_p(t) \right) \cdot 
		\frac{d \mathbf{M}_c^{-1}}{dt} +
		\left( 
		\sum \limits_{p \in \mathcal{P}(c)}  \, \frac{\partial \psi_p}{\partial \boldsymbol{\zeta} }
		\, 
		\frac{d Z_p}{dt} \right)
		\cdot \mathbf{M}_c^{-1}
		\nonumber \\
		&=& 
		-\left( 
		\sum \limits_{p \in \mathcal{P}(c)}  \, \frac{\partial \psi_p}{\partial \boldsymbol{\zeta} }
		\, 
		Z_p(t) \cdot \mathbf{M}_c^{-1} \right) \cdot 
		\left(  \sum \limits_{p \in \mathcal{P}(c)}  \, \frac{\partial \psi_p}{\partial \boldsymbol{\zeta} }
		\, 	\vv_p \, \mathbf{M}_c^{-1} \right) 
		\nonumber \\ 
		&& -
		\left( 
		\sum \limits_{p \in \mathcal{P}(c)}  \, \frac{\partial \psi_p}{\partial \boldsymbol{\zeta} }
		\, 
		\auxf_p \right)
		\cdot \mathbf{M}_c^{-1}
		= 
		-\nabla_c^p Z_p \cdot \nabla_c^q \vv_q - \nabla_c^p \auxf_p,
		\label{eqn.dJdt}   
	\end{eqnarray} 
	which is identical with \eqref{eqn.J.compact} and thus with \eqref{eqn.sdhtc_J} since $\J_c(0) = \nabla_c^p Z_p(0))$ holds. As such we can conclude that $\J_c = \nabla_c^p Z_p$. Application of the 
	discrete curl operator and assuming exact time integration leads to the sought result 
	\begin{equation}
		  \nabla_p^c \times \J_c = \nabla_p^c \times \nabla_c^q Z_q = 0. 
	\end{equation}
\end{proof}

\section{Numerical results} \label{sec.results}
We present a set of test cases that show the accuracy of the new Lagrangian HTC scheme for the GPR model, which will be labeled as Lag-HTC. Inviscid and viscous fluids are considered as well as ideal elastic solids. For comparison, the Lagrangian scheme \cite{LGPR} is used to run some of the test cases, which is an extension of the EUCCLHYD scheme \cite{EUCCLHYD} to the GPR model. The detailed equations of state (EOS) and the functions $\theta_1(\tau_1)$ and $\theta_2(\tau_2)$ used in this paper are given in \cite{GPRmodel,SIGPR}.
The polytropic index of the gas is assumed to be $\gamma=7/5$ and the specific heat at constant volume is $c_v=2.5$. Whenever a viscosity coefficient $\mu$ is specified, the relaxation time $\tau_1$ is computed according to $\mu = \frac{1}{6} \rho_0 c_s^2 \tau_1$. Likewise, if a heat conduction coefficient $\kappa \neq 0$ is prescribed, the corresponding relaxation time $\tau_2$ is evaluated from the asymptotic relation $\kappa = \rho_0 T_0 c_h^2 \tau_2$. In the other cases, no source terms are considered, thus we set $\tau_1=\tau_2=10^{20}$. The distortion matrix is always initialized as $\A=\mathbf{I}$ since $\x(0)=\boldsymbol{\xi}$, and the thermal impulse is initially $\J=\mathbf{0}$. The reference density is given by $\rho_c(0)=\rho(\x_c,t=0)$. If not stated otherwise, we set $\text{CFL}=0.05$ and a six-stage fifth order Runge-Kutta scheme is adopted for time integration, thus ensuring that the GCL and total energy conservation are satisfied up to sufficient accuracy at the fully-discrete level.
In the numerical simulations, we monitor the following properties over time.
\begin{itemize}[leftmargin=+0.3in]
	\item Curl operators: the error related to the curl of the thermal impulse vector $\J$ is computed at a given time $t=t^n$ by considering 
	\begin{equation*}
	\epsilon^{\nabla \times \J} := \max \limits_{p} \left|	- \frac{1}{|\omega_p|} \sum \limits_{c \in \mathcal{C}(p)} \lnpc \times \J_c \right|.
	\end{equation*} 
    The error related to the curl of the distortion matrix $\A$ is computed by applying the above relation to each row of $\A$, and then taking the maximum norm.
	\item Determinant of $\A$: the error related to the geometric constraint $|\A_c|=\rho_c(t)/\rho_c(0)$ is computed at each time step as
	\begin{equation*}
		  \epsilon^{\detA} := \max \limits_c \left|  -\frac{|\A_c|}{|\omega_c|} \sum \limits_{p \in \mathcal{P}(c)} \lnpc \cdot ({\vv}_p-\vv_c) + \frac{\rho_c^2}{\rho_c(0)} \frac{d \tau_c}{dt} \right|.
	\end{equation*}	
	\item Total energy conservation: the error at each time step is given by
	\begin{equation*}
		\epsilon^{E}:= \left| \sum_{c} \w_c \cdot \mde{\q}{t} \right|.
	\end{equation*}
\end{itemize}

\subsection{Numerical convergence study}
\begin{figure}[!htbp]
	\begin{center}
		\begin{tabular}{cc}
			\includegraphics[trim= 5 5 5 5, clip,width=0.47\textwidth]{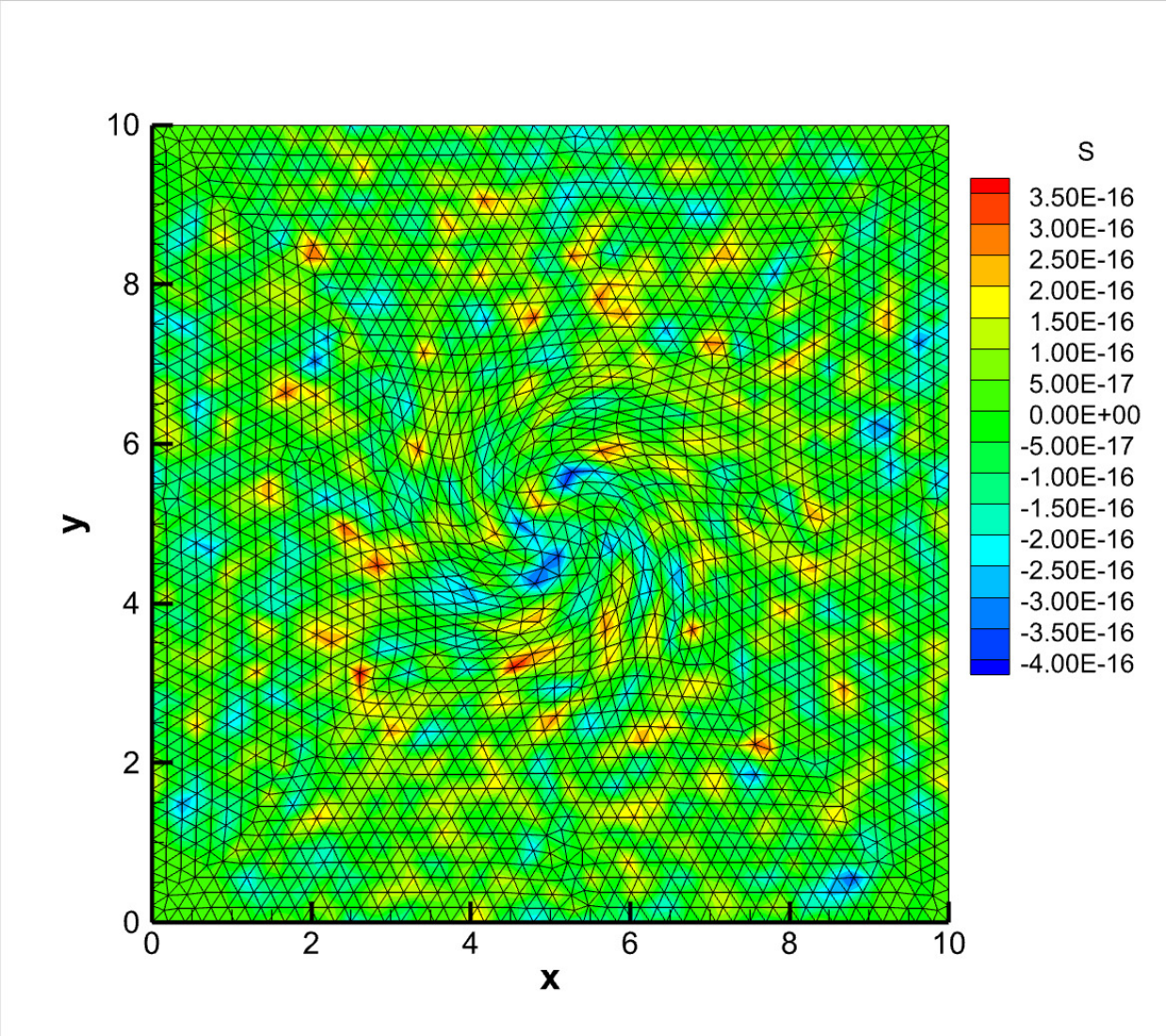}  &
			\includegraphics[trim= 5 5 5 5, clip,width=0.47\textwidth]{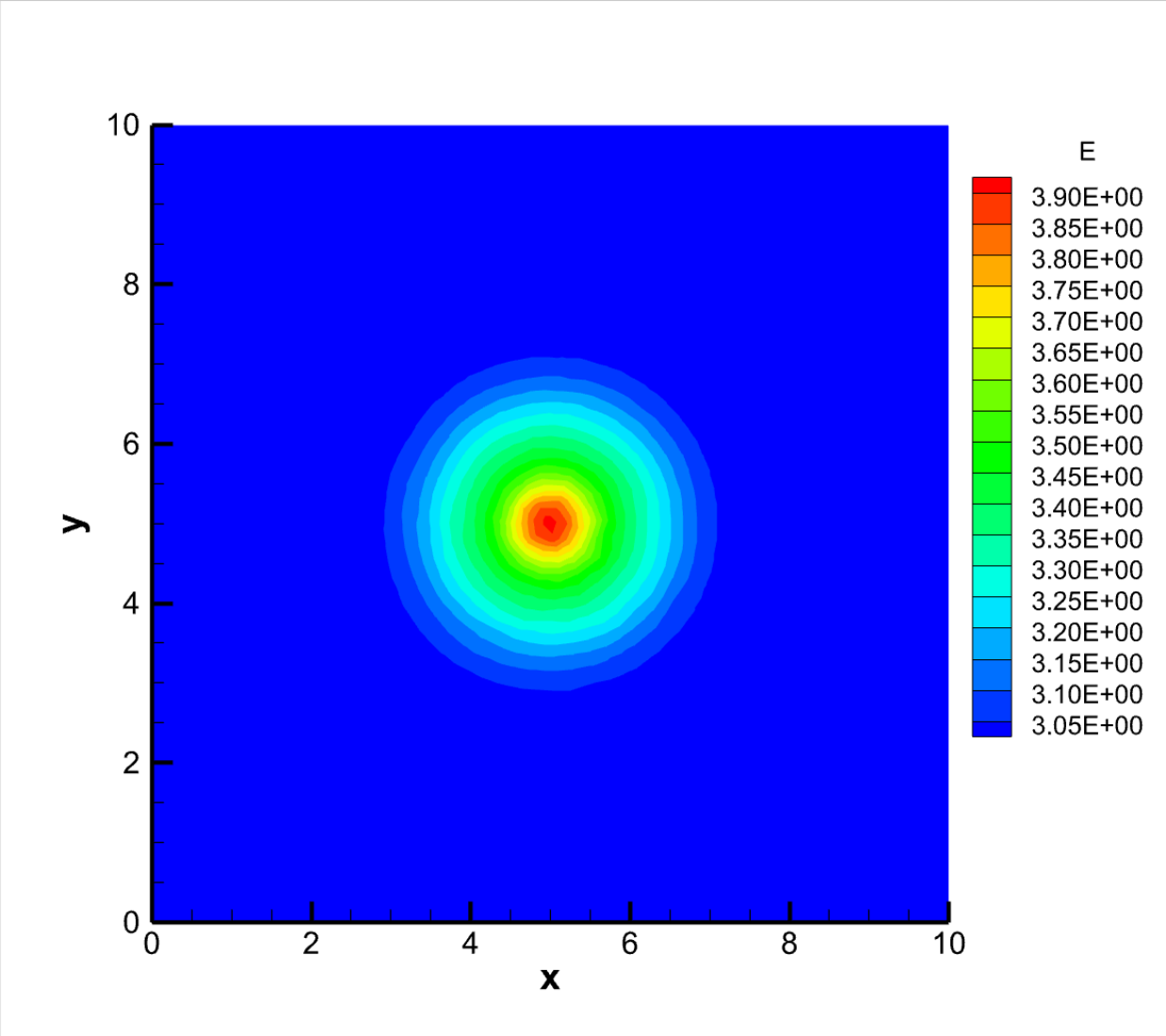} \\
		\end{tabular}
		\caption{Isentropic vortex problem at time $t_f=1$. Left: entropy distribution and mesh configuration. Right: energy distribution.}
		\label{fig.ShuVortex}
	\end{center}
\end{figure}

\begin{figure}[!htbp]
	\begin{center}
		\begin{tabular}{cc}
			\includegraphics[trim= 5 5 5 5, clip,width=0.47\textwidth]{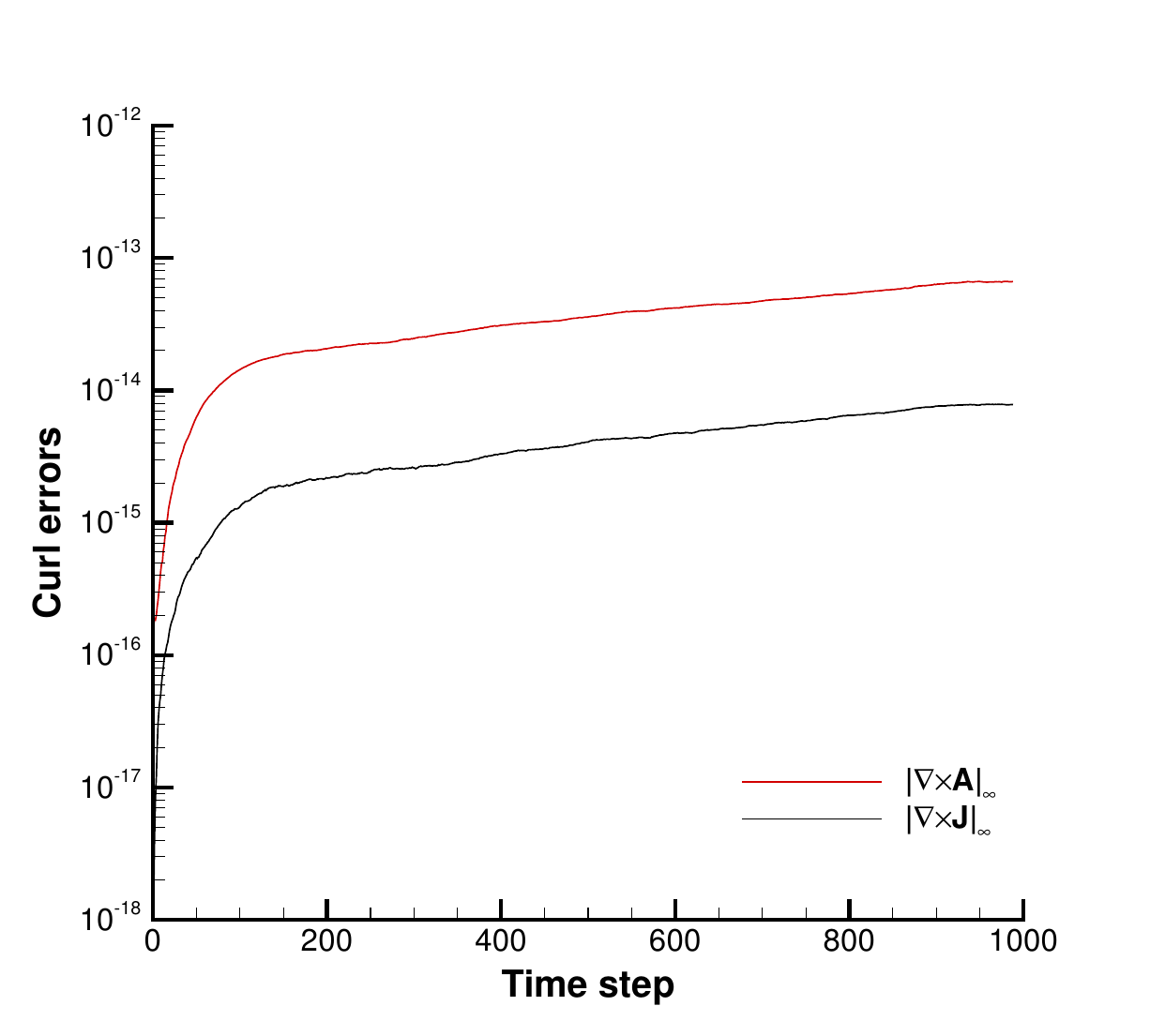}  &
			\includegraphics[trim= 5 5 5 5, clip,width=0.47\textwidth]{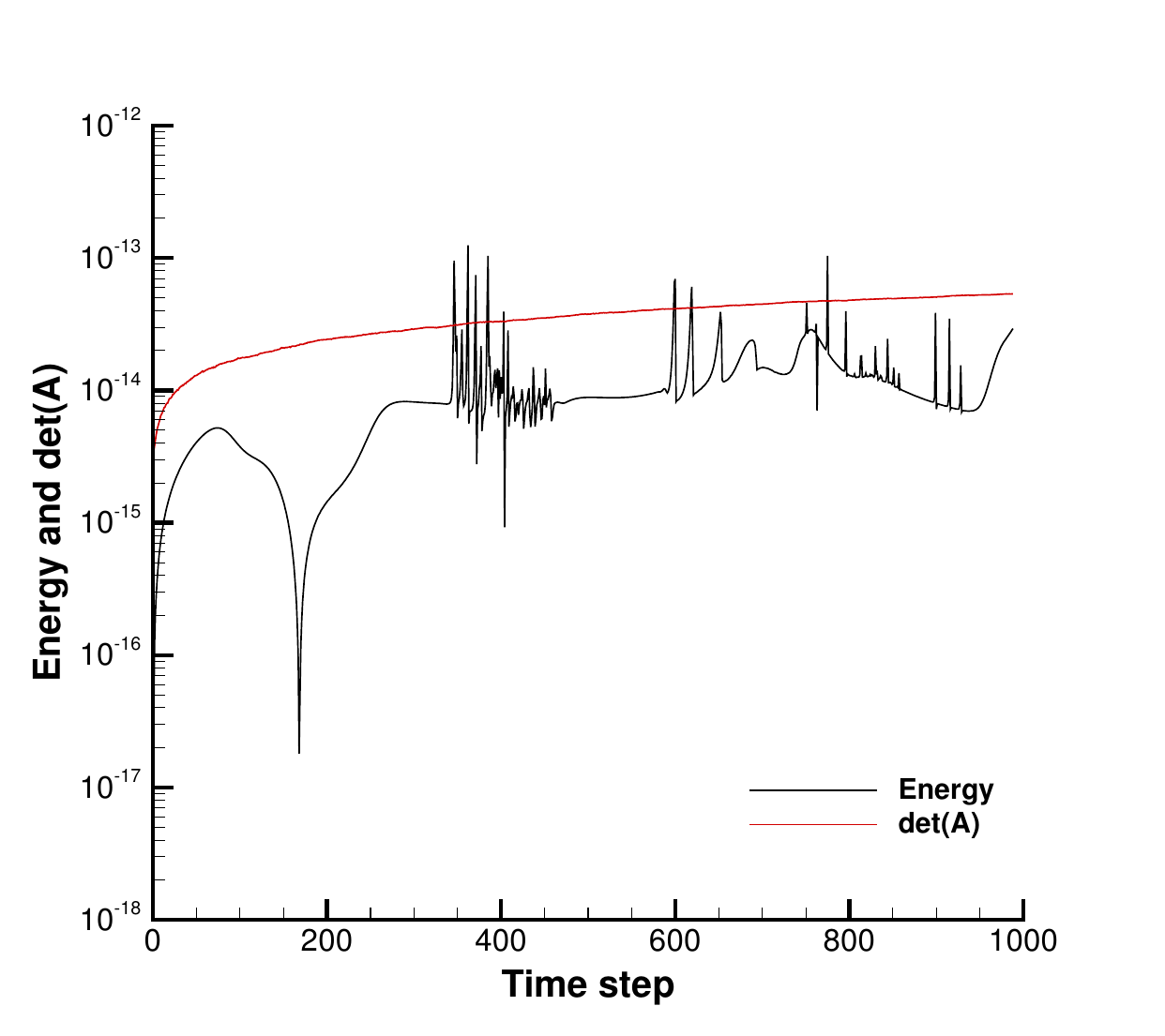} \\
		\end{tabular}
		\caption{Isentropic vortex problem. Left: time evolution of the curl of $\A$ and $\J$. Right: monitoring of the determinant of $\A$ compatibility and total energy conservation error.}
		\label{fig.ShuVortex_curl-detA}
	\end{center}
\end{figure}

We solve the isentropic vortex problem \cite{HuShuTri} to assess the accuracy of the scheme. We set $c_s=0$ and $c_h=0$ to model the inviscid Euler equations of compressible gas dynamics. The computational domain is $\Omega(0)=[0;10]^2$ with slip wall boundaries everywhere. The initial condition is detailed in \cite{HuShuTri}, and here the vortex is not transported, thus the convective velocity is set to zero. This test case is run until the final time $t_f=1$ on a sequence of successively refined meshes of characteristic mesh size 
$ 
	h = \max \limits_{\Omega(t_f)} \frac{|\omega_c|}{|\partial \omega_c|}.
$
The experimental order of convergence (EOC) is reported in Table \ref{tab.convRates}, showing that the expected first order of accuracy is retrieved. The errors are measured in $L_2$ norm for density, horizontal velocity and pressure. Figure \ref{fig.ShuVortex} shows the entropy and the energy distribution at time $t_f=1$ on a computational mesh composed of $N_c=5180$ cells. As expected, the flow remains exactly isentropic, while conserving at the same time total energy. The time evolution of the curl of $\A$ and $\J$ are depicted in Figure \ref{fig.ShuVortex_curl-detA}, where we also monitor the error related to the determinant of $\A$ and the total energy conservation. The curl-free and the algebraic constraints are well-preserved by the scheme. 

\begin{table}[!htbp]  
	\caption{Experimental order of convergence (EOC) for the isentropic vortex problem using the Lag-HTC scheme. The errors are measured in the $L_2$ norm and refer to the variables $\rho$ (density), $u$ (horizontal velocity) and $p$ (pressure) at time $t_{f}=1$.}  
	\begin{center} 
		\begin{small}
			\renewcommand{\arraystretch}{1.2}
			\begin{tabular}{c|cccccc}
				$h$ & $L_2(\rho)$ & EOC($\rho$) & $L_2(u)$ & EOC($u$) & $L_2(p)$ & EOC($p$) \\ 
				\hline
				3.254E-01 & 2.6883E-01 & -    & 1.4736E-01 & -    & 3.6822E-01 & -    \\
				2.490E-01 & 2.1588E-01 & 0.82 & 1.1103E-01 & 1.06 & 2.9549E-01 & 0.82 \\
				1.654E-01 & 1.4773E-01 & 0.93 & 7.3479E-02 & 1.01 & 2.0149E-01 & 0.94 \\
				1.283E-01 & 1.1088E-01 & 1.13 & 5.6521E-02 & 1.03 & 1.5117E-01 & 1.13 \\
			\end{tabular}
		\end{small}
	\end{center}
	\label{tab.convRates}
\end{table} 

\subsection{Viscous shock wave}
\begin{figure}[!htbp]
	\begin{center}
		\begin{tabular}{cc}
			\includegraphics[trim= 5 5 5 5, clip,width=0.47\textwidth]{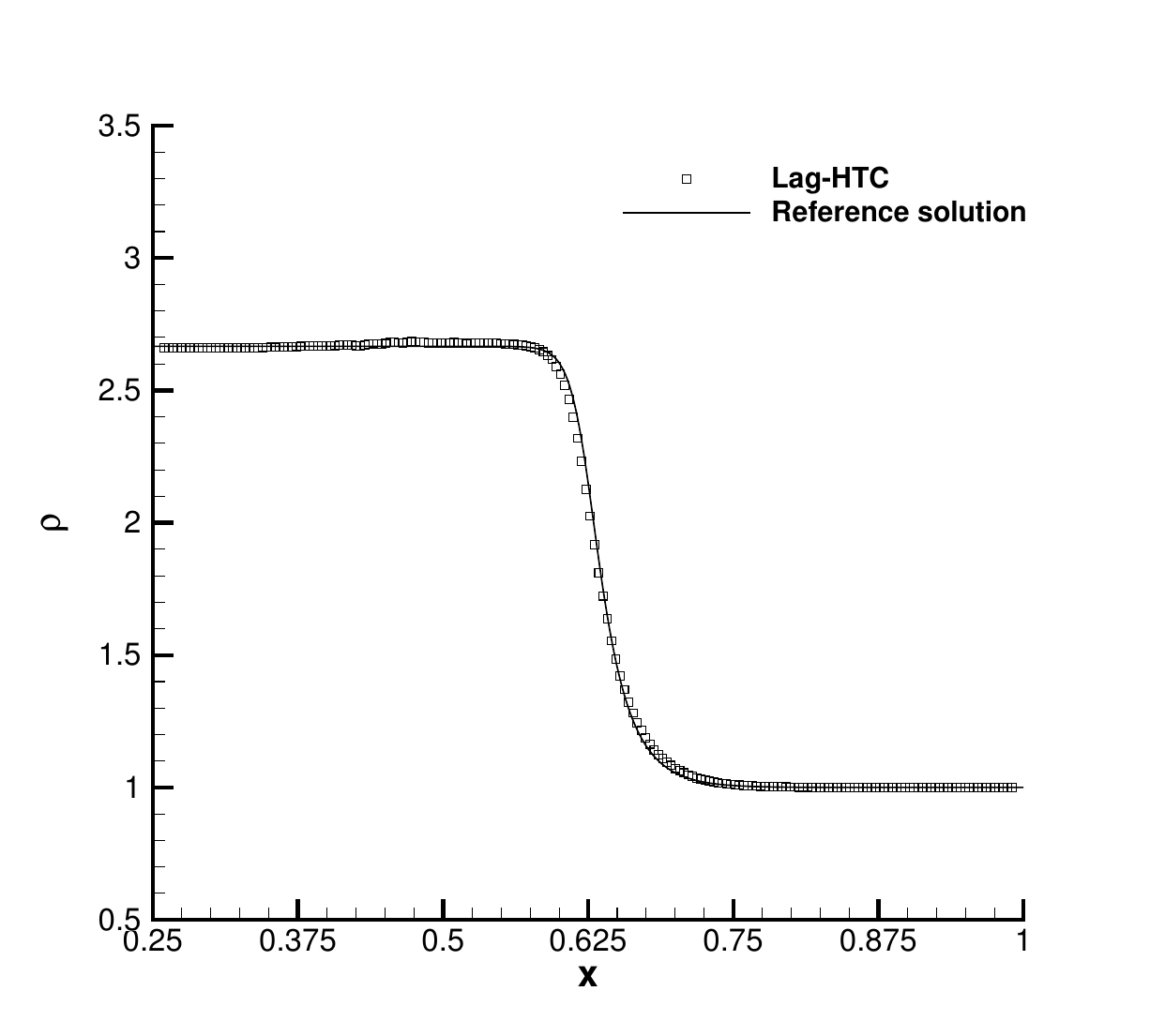}  &
			\includegraphics[trim= 5 5 5 5, clip,width=0.47\textwidth]{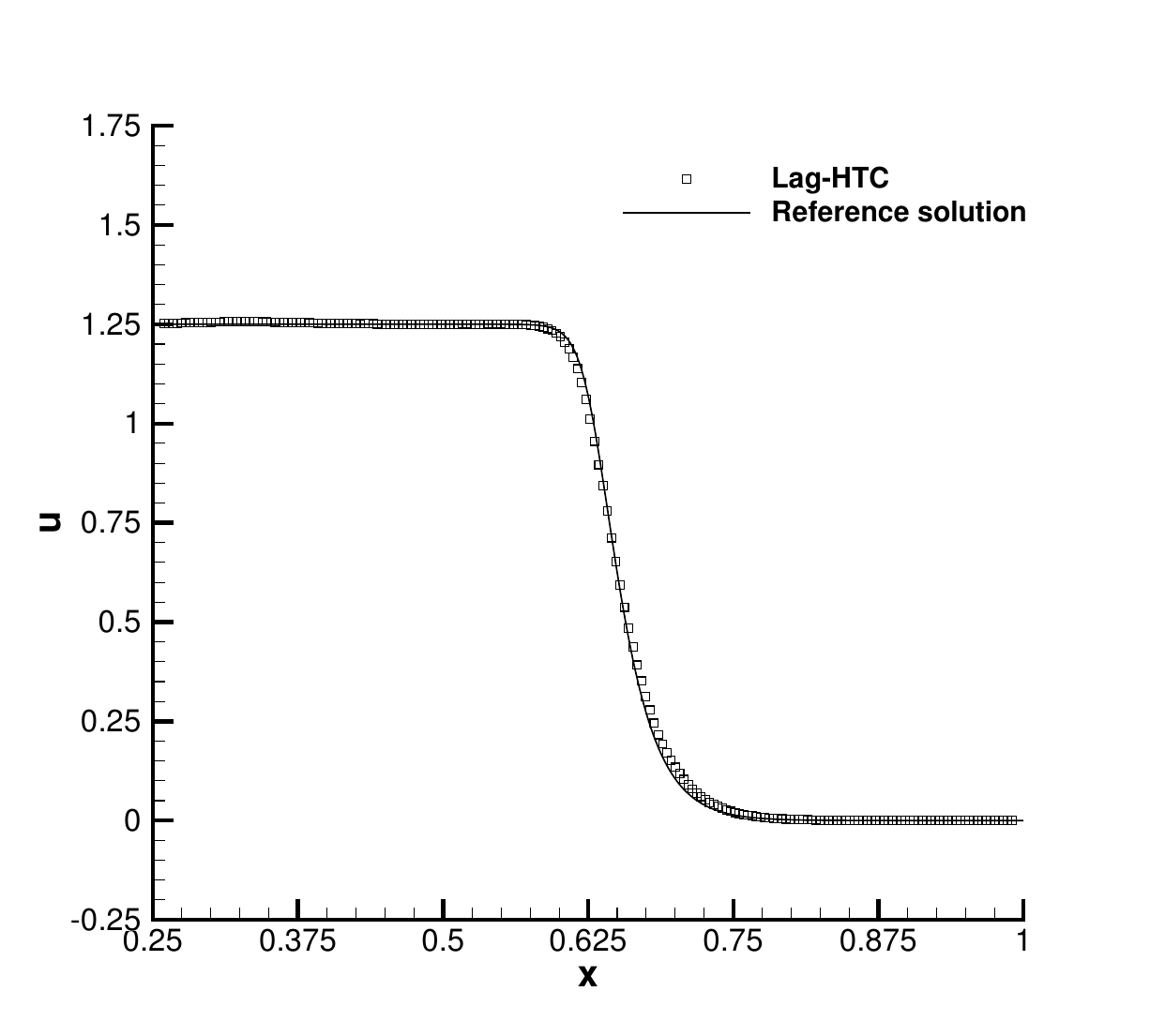} 
		\end{tabular}
		\caption{Viscous shock problem at time $t_f=0.2$. Numerical results for density (left) and horizontal velocity (right) 
			compared against the exact solution.}
		\label{fig.viscousshock}
	\end{center}
\end{figure}

\begin{figure}[!htbp]
	\begin{center}
		\begin{tabular}{cc}
			\includegraphics[trim= 5 5 5 5, clip,width=0.5\textwidth]{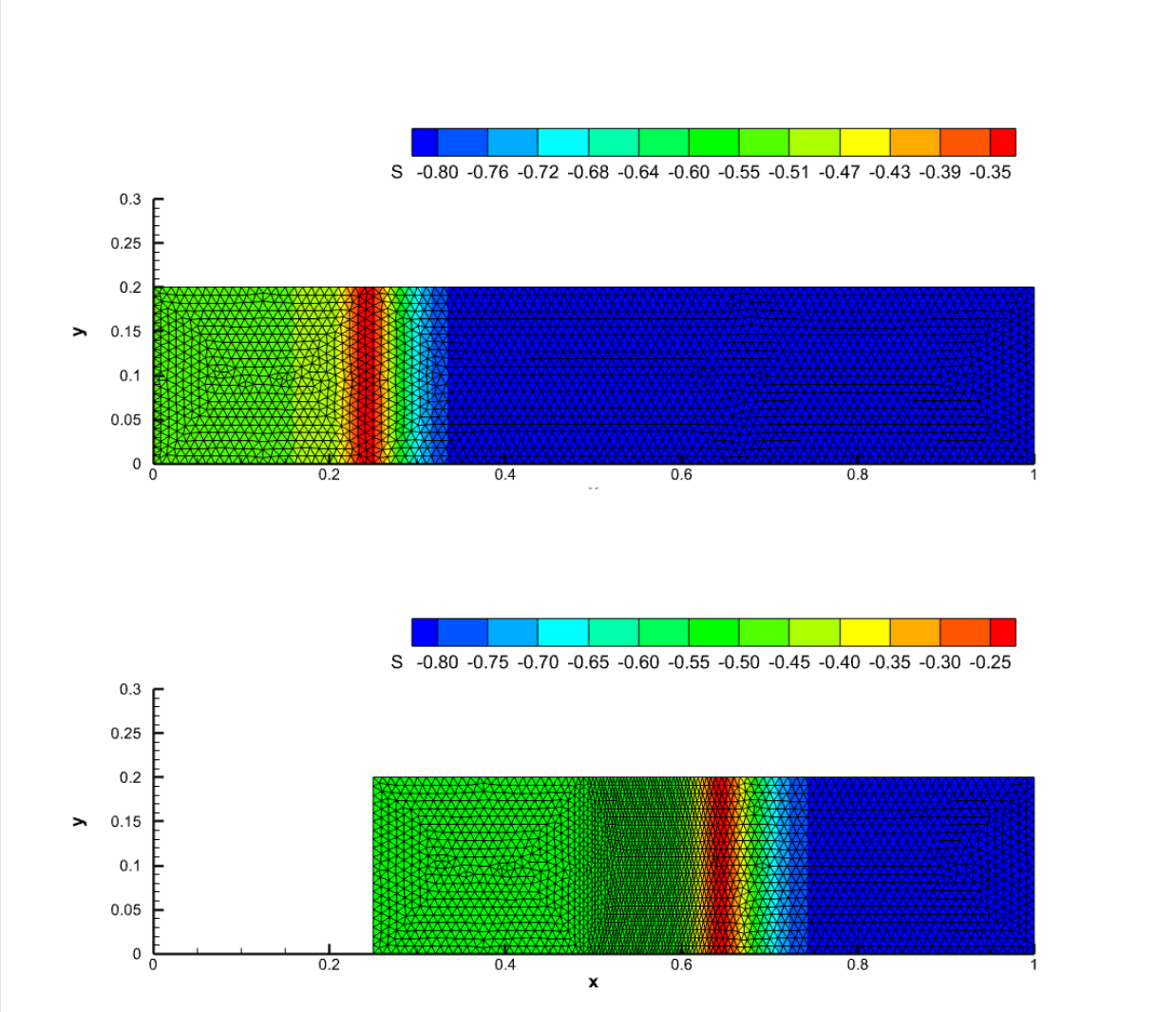}  &
			\includegraphics[trim= 5 5 5 5, clip,width=0.47\textwidth]{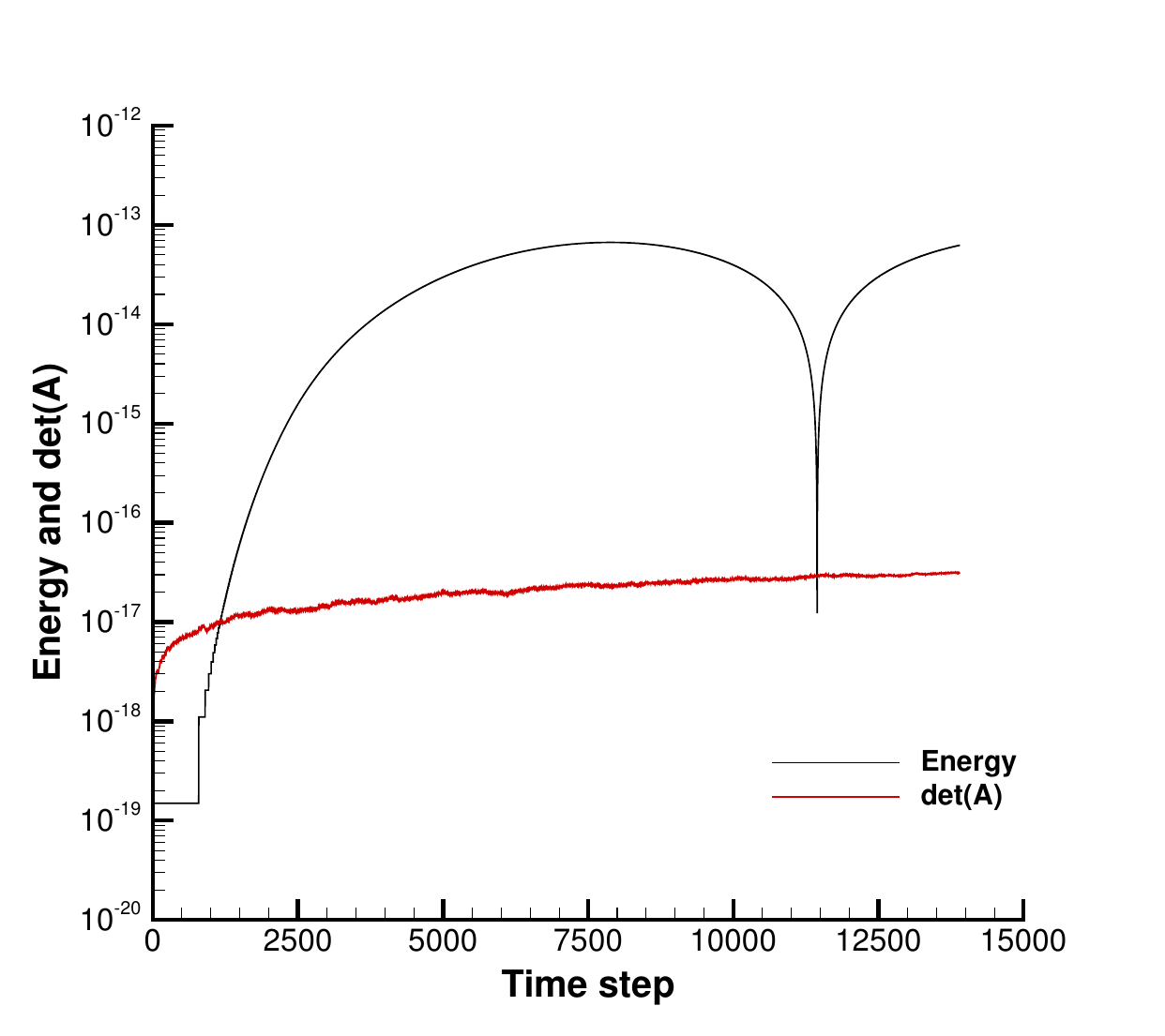} \\
		\end{tabular}
		\caption{Viscous shock problem. Left: mesh configuration and entropy distribution at the initial time (top) and at the final time (bottom). Right: monitoring of the determinant of $\A$ compatibility and total energy conservation error.}
		\label{fig.viscousshock_S-detA}
	\end{center}
\end{figure}

We model compressible heat-conducting viscous flows by setting $c_s=c_h=20$, $\mu=2 \cdot 10^{-2}$ and $\kappa=9.3333 \cdot 10^{-2}$. The initial computational domain is  $\Omega(0)=[0;1] \times [0;0.2]$, discretized with triangles of size $h=1/100$. No-slip wall boundaries are imposed in $y-$direction, while a Dirichlet boundaries are prescribed in $x$ direction. For Prandtl number $\text{Pr}=0.75$, there exists an exact solution, see \cite{Becker1923}. The Reynolds number is $\text{Re}_s=\rho_0 c_0 M_s L \mu^{-1}$, with the reference length $L=1$. The initial condition is a shock centered at $x=0.25$, propagating at Mach $M_s=2$ with $\text{Re}_s=100$. The upstream shock state is defined by 
$
	(\rho,u,v,p) = (\rho_0, \, 0, \, 0, \, 1/\gamma ), 
$ 
with $c_0=1$ and $\rho_0=1$. 
The comparison between the numerical solution obtained with the Lagrangian HTC scheme and the exact solution of the compressible Navier-Stokes equations at time $t_f=0.2$ is shown in \ref{fig.viscousshock}. For all quantities a very good agreement is achieved. The initial and final mesh configuration with the entropy distribution is depicted in Figure \ref{fig.viscousshock_S-detA}, where we see that entropy is increasing while the shock is moving into the fluid. Furthermore, we also plot the error related to the determinant of $\A$ and the total energy conservation. Both constraints are preserved up to machine precision.

\subsection{Solid rotor problem}
\begin{figure}[!htbp]
	\begin{center}
		\begin{tabular}{cc}
			\includegraphics[trim= 5 5 5 5, clip,width=0.47\textwidth]{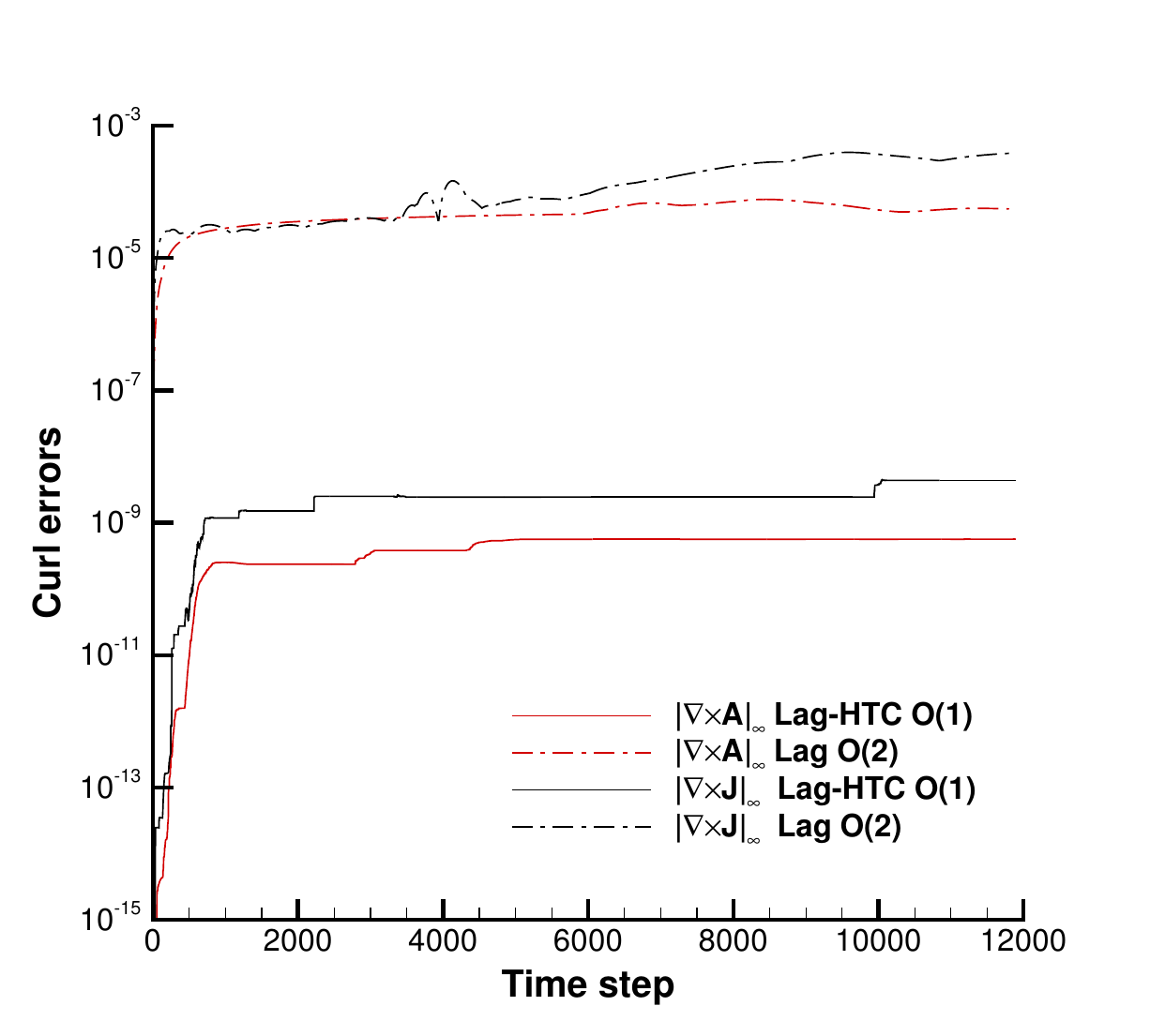}  &
			\includegraphics[trim= 5 5 5 5, clip,width=0.47\textwidth]{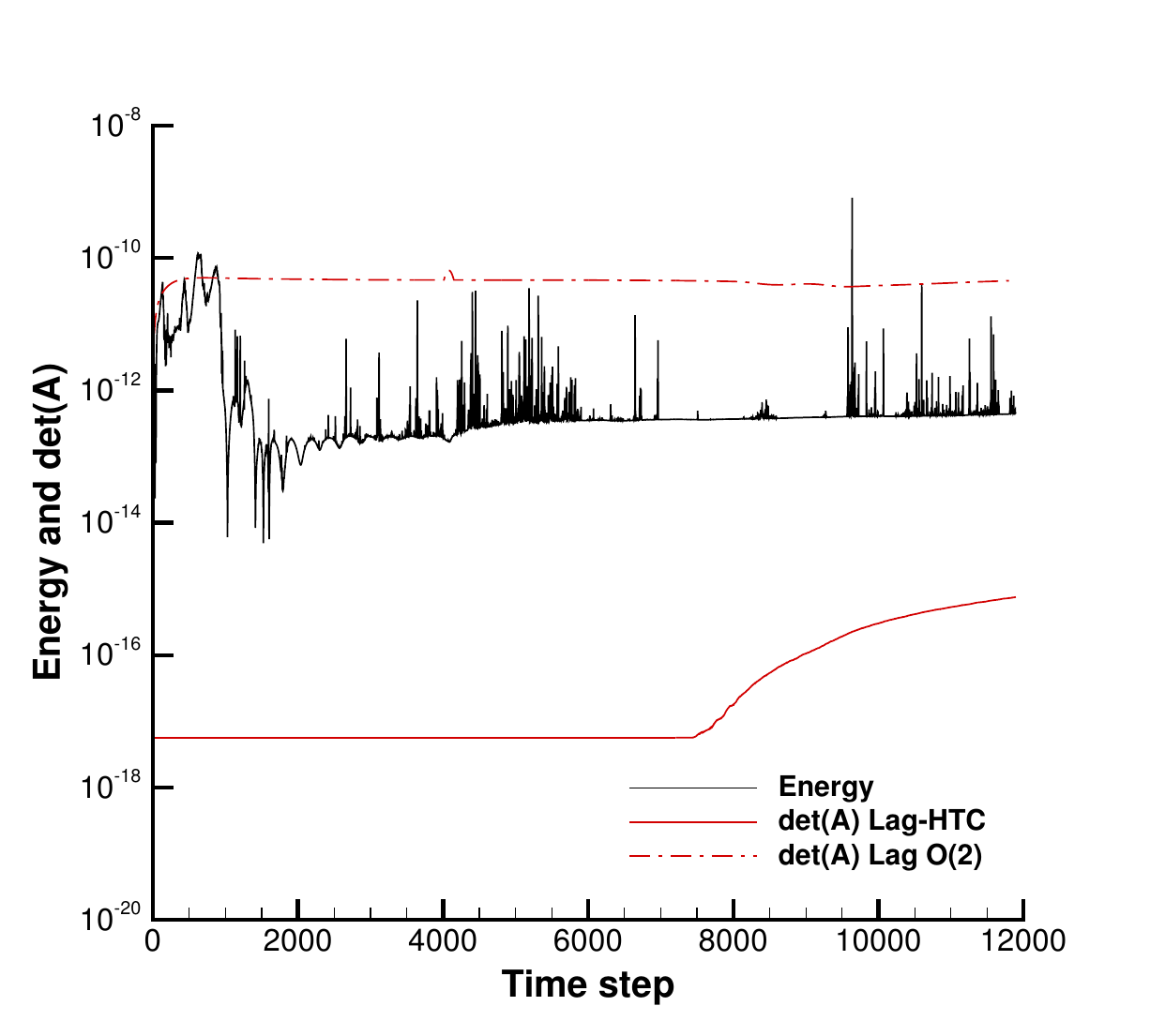} \\
		\end{tabular}
		\caption{Solid rotor problem. Left: time evolution of the curl of $\A$ and $\J$. Right: monitoring of the determinant of $\A$ compatibility and total energy conservation error. Comparison with the second order EUCCLHYD scheme of \cite{LGPR}.}
		\label{fig.SolidRotor_curl-detA}
	\end{center}
\end{figure}

The solid rotor problem \cite{SIGPR,HTCGPR} deals with ideal elastic solids. We fix $c_s=c_h=1$ and the final time of the simulation is $t_f=0.3$. The initial computational domain is $\Omega(0)=[-1;1]^2$ with slip wall boundaries, which is discretized with a characteristic mesh size of $h=1/150$ yielding a total number of $N_c=94754$ cells. The initial condition of the material reads
\begin{equation}
	(\rho,u,v,p)=\left\{
	\begin{array}{lll}
		(1, \, -y/R, \, x/R, \, 1) & & r<R \\
		(1,\, 0,\, 0,\, \, 1) & &  r\geq R
	\end{array} \right., \qquad t=0, \quad \x \in \Omega(0),
\end{equation}
with the initial discontinuity located at $R=0.2$ and $r=\sqrt{x^2+y^2}$. 
%
Figure \ref{fig.SolidRotor_curl-detA} shows the time evolution of the curl of $\A$ and $\J$ as well as the error related to the determinant of $\A$ and the total energy conservation over time. We can conclude that all the properties of the semi-discrete Lag-HTC scheme are fulfilled also for this test case.

\subsection{Elastic vibration of a beryllium plate}
\begin{figure}[!htbp]
	\begin{center}
		\begin{tabular}{cc}
			\includegraphics[trim= 5 5 5 5, clip,width=0.47\textwidth]{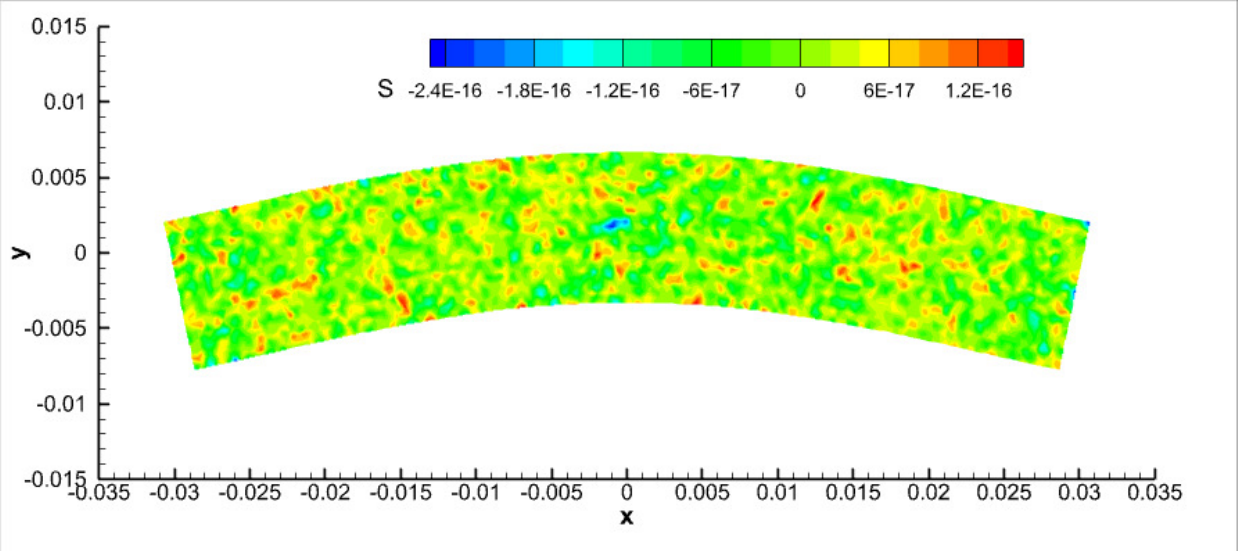}  &
			\includegraphics[trim= 5 5 5 5, clip,width=0.47\textwidth]{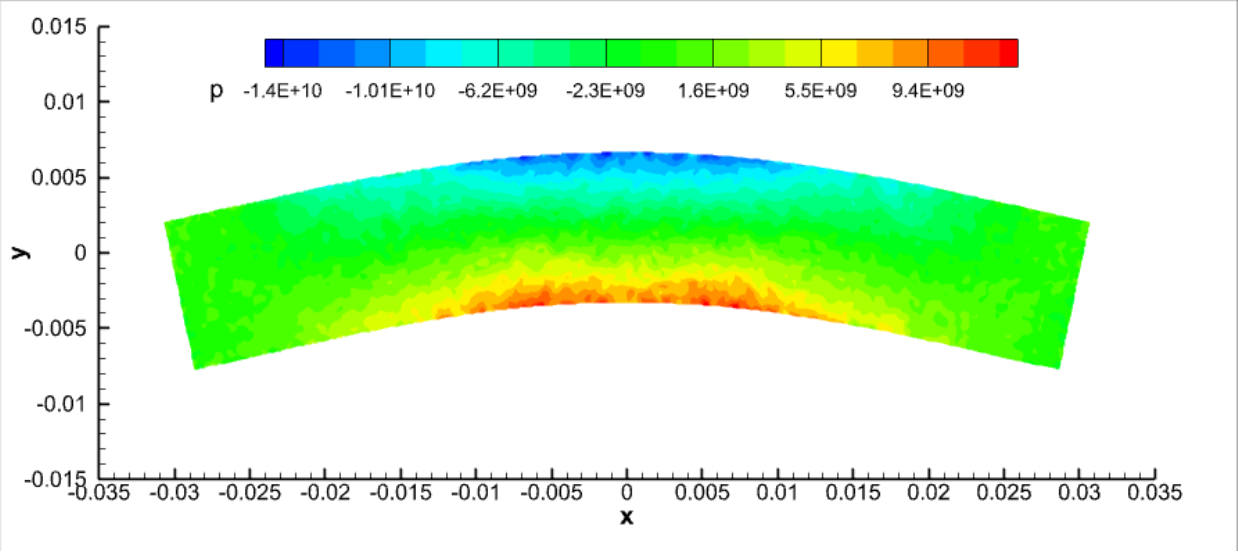} \\
			\includegraphics[trim= 5 5 5 5, clip,width=0.47\textwidth]{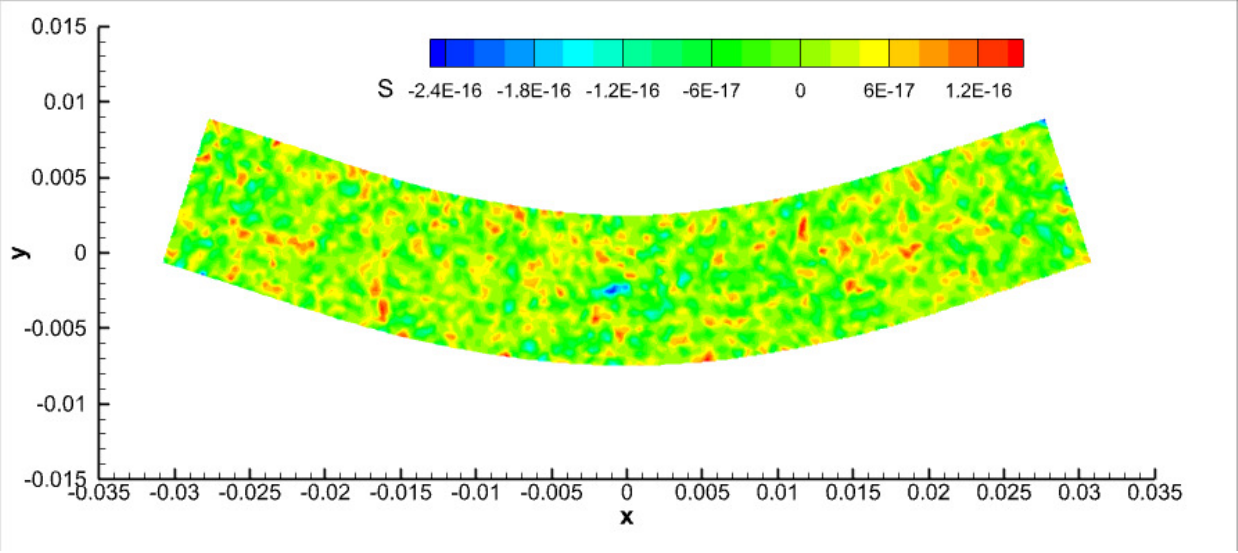}  &
			\includegraphics[trim= 5 5 5 5, clip,width=0.47\textwidth]{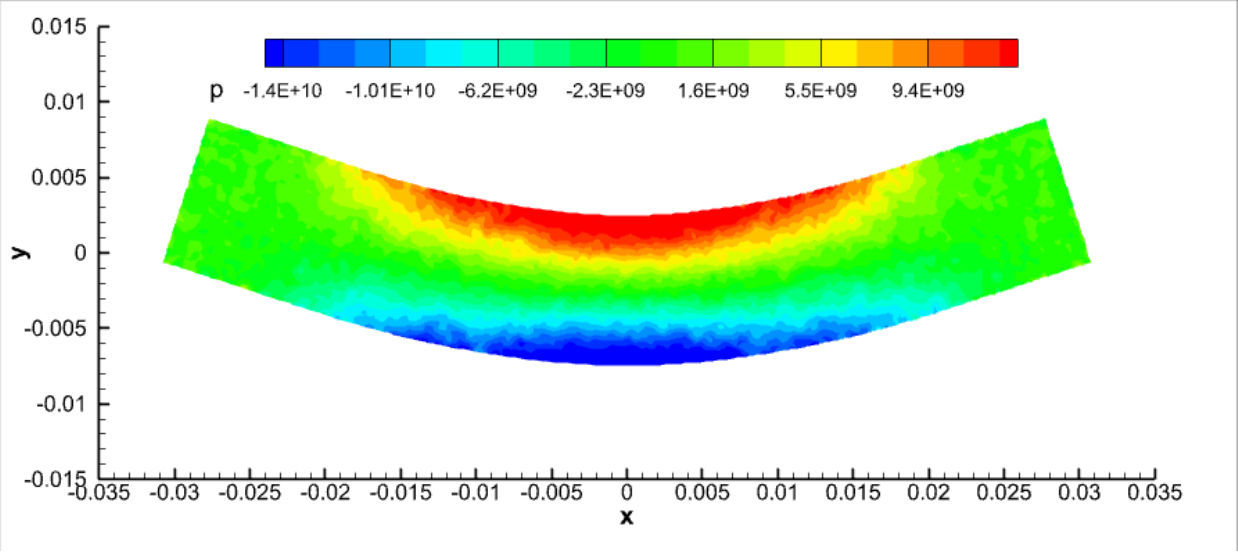} \\
		\end{tabular}
		\caption{Beryllium plate problem. Entropy (left) and pressure (right) distribution at time $t=1\cdot 10^{-5}$ (top) and $t=2\cdot 10^{-5}$ (bottom).}
		\label{fig.BePlate_S-p}
	\end{center}
\end{figure}

\begin{figure}[!htbp]
	\begin{center}
		\begin{tabular}{cc}
			\includegraphics[trim= 5 5 5 5, clip,width=0.47\textwidth]{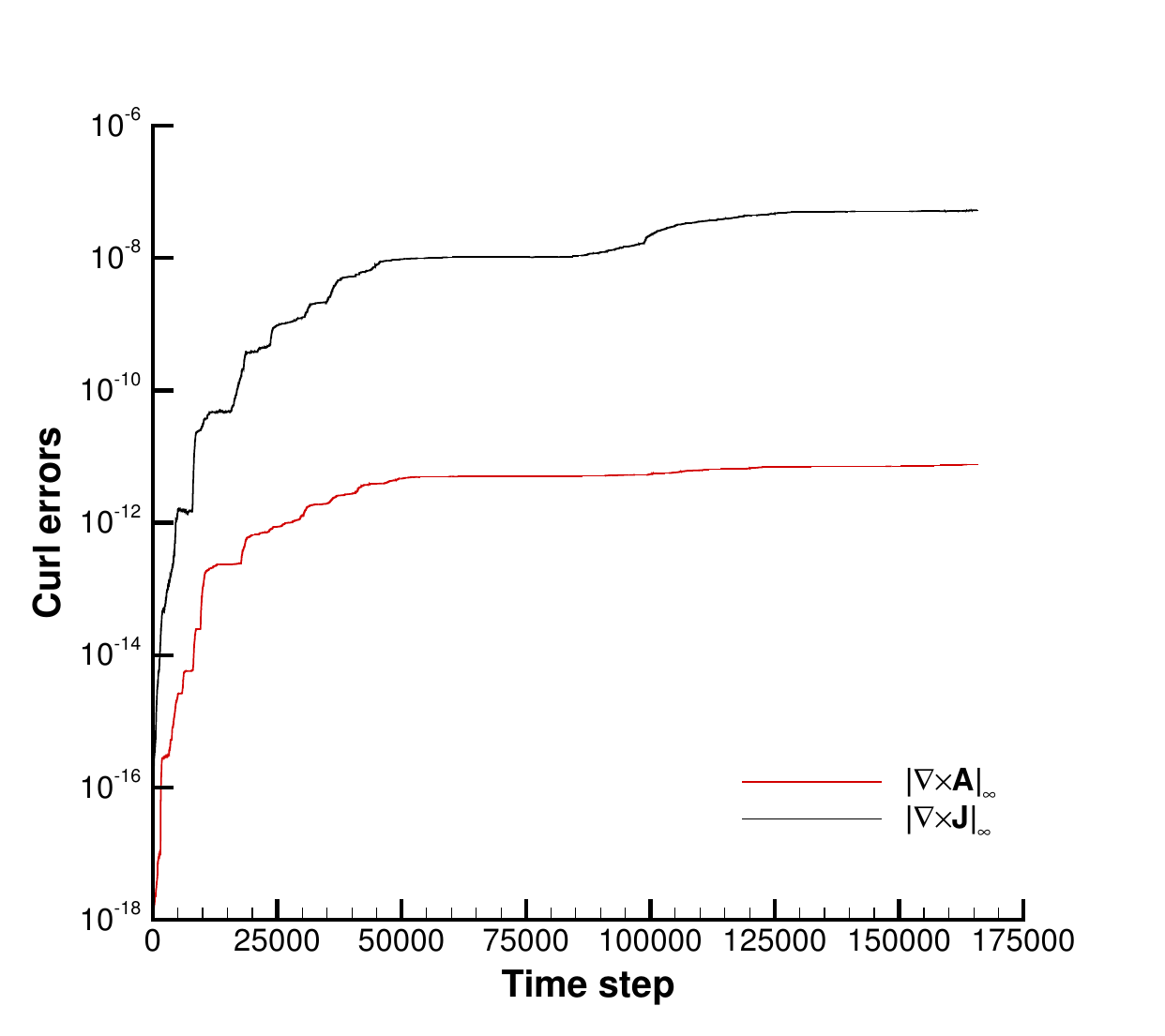}  &
			\includegraphics[trim= 5 5 5 5, clip,width=0.47\textwidth]{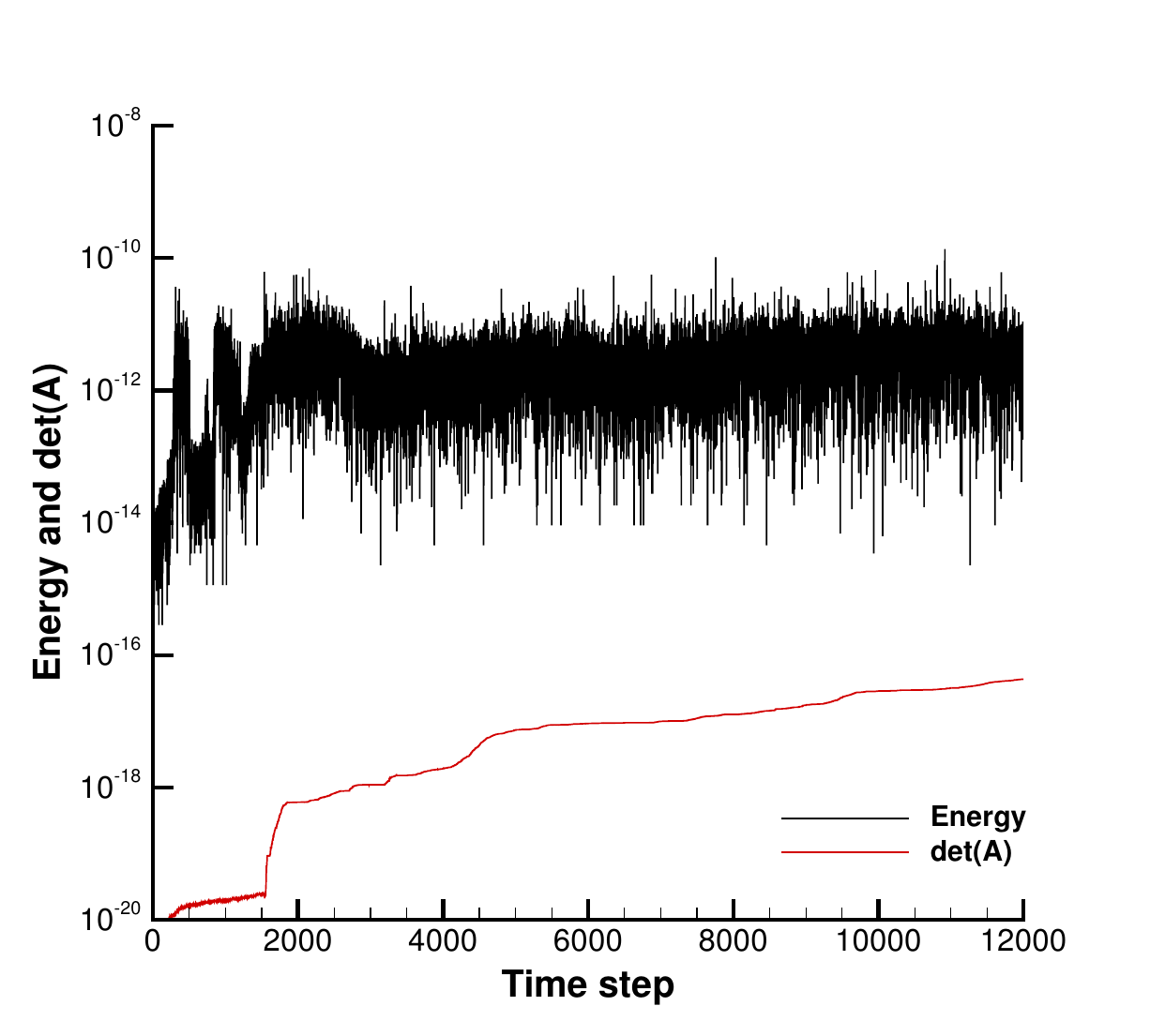} \\
		\end{tabular}
		\caption{Beryllium plate problem. Left: time evolution of the curl of $\A$ and $\J$. Right: monitoring of the determinant of $\A$ compatibility and total energy conservation error.}
		\label{fig.BePlate_curl-detA}
	\end{center}
\end{figure}

\begin{figure}[!htbp]
	\begin{center}
		\begin{tabular}{c}
			\includegraphics[trim= 10 60 10 10, clip,width=0.6\textwidth]{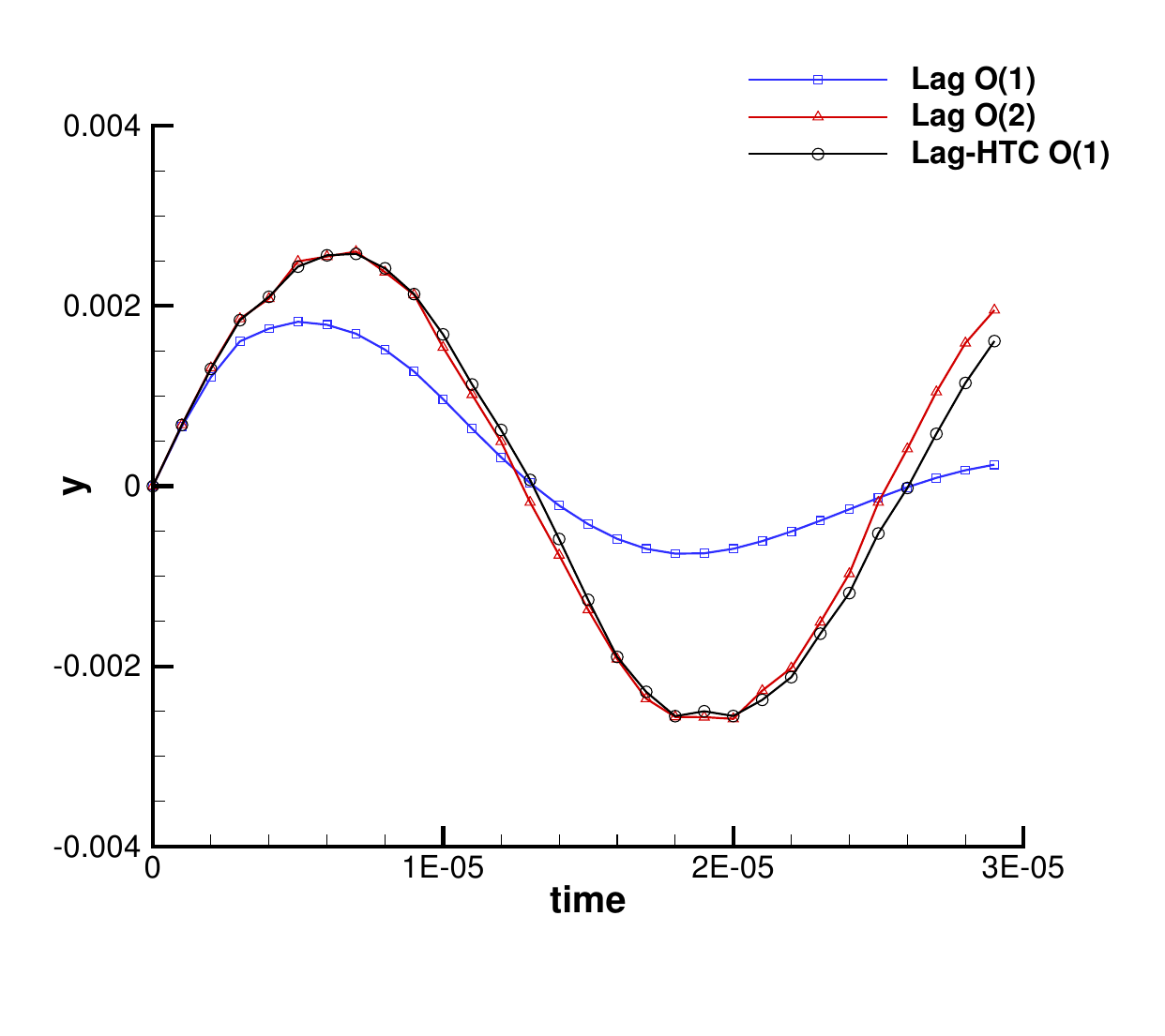} \\
		\end{tabular}
		\caption{Beryllium plate problem. Vertical displacement of the point initially located at $\x = (0, 0)$ with first and second order EUCCLHYD scheme \cite{LGPR} as well as with the new Lag-HTC scheme.}
		\label{fig.BePlate_disp}
	\end{center}
\end{figure}

This test case describes the elastic vibration of a beryllium plate, see \cite{Kluth10,ShashkovCellCentered}. The 
 domain is $\Omega(0)=[-0.03;0.03]\times[-0.005;0.005]$ and the length of the plate is $L=0.06$. The mesh is composed of $N_c=8414$ elements, and free-traction boundary conditions are set everywhere. The final time is chosen as $t_f=3\cdot 10^{-5}$, so that approximately one oscillating period is completed. We use a stiffened gas equation of state with $\gamma=1.11$ and $\Gamma=1.124$. The other parameters of the GPR model are: $c_s=9046.59$, $c_0=12870$, $c_h=0$, $c_v=1$. The material initially has uniform density $\rho(0)=1$ and pressure $p(0)=0$, and is perturbed via a vertical velocity field of the form
\begin{equation}
	v(x) = A \omega \left[  a_1(\sinh(x')+\sin(x')) - a_2(\cosh(x')+\cos(x')) \right],
\end{equation}
where $x'=\alpha(x+L/2)$, $\alpha=78.834$, $A=4.3369\times 10^{-5}$,
$\omega=2.3597\times 10^5$, $a_1=56.6368$ and $a_2=57.6455$. We show the entropy and pressure distribution in Figure \ref{fig.BePlate_S-p} at different times. One can notice that the purely elastic behavior of the solid does not produce any entropy within our new Lagrangian HTC scheme. Figure \ref{fig.BePlate_curl-detA} confirms the preservation of the curl-free property as well as the geometric relation on $\detA$ and total energy conservation over more than $170'000$ time steps. Finally, this test problem is run with the first and second order EUCCLHYD scheme devised in \cite{LGPR}. In Figure \ref{fig.BePlate_disp} we show the time evolution of the vertical displacement of the point originally located at $\x = (0, 0)$. The first and second order EUCCLHYD schemes are compared against the new first order Lagrangian HTC scheme, showing that the thermodynamically compatible method substantially reduces numerical dissipation and the results are comparable to a second order EUCCLHYD scheme. 

\section{Conclusions} \label{sec.concl}

In this paper we have presented a new numerical method for the solution of the unified GPR model of continuum mechanics. The scheme preserves all essential structural properties of the governing PDE system exactly also on the discrete level. More precisely, the method is by construction compatible with the second principle of thermodynamics, since the entropy inequality is directly discretized within our scheme. Furthermore, thanks to a judicious choice of the nodal fluxes, the method also conserves total energy at the semi-discrete level and thus is also compatible with the first principle of thermodynamics. The fact that we directly discretize the entropy inequality and obtain the energy conservation law as a consequence is in line with the SHTC formalism of Godunov and Romenski \textit{et al.}, but is still rather unconventional at the discrete level. In addition, our new schemes are also by construction exactly compatible with the determinant constraint on the distorsion field. This algebraic geometric identity is not easy to achieve at the discrete level. In our scheme, the discrete compatibility is possible thanks to a judicious choice of nodal fluxes and compatible nabla operators. Last but not least, in the absence of algebraic source terms, the method presented in this paper also preserves the curl-free property of the distorsion field $\A$ and of the thermal impulse $\J$ exactly at the semi-discrete level. As such, all basic structural properties of the PDE system are satisfied by our discretization. 
The numerical results obtained for some basic test cases underline that the properties are also achieved in practice. 
Future work will concern the extension to higher order of accuracy, following the approach recently outlined in \cite{CompatibleDG1}. 

\section*{Acknowledgments}
This research was funded by the Italian Ministry of University and Research (MUR) via the PRIN Project 2022 No. 2022N9BM3N,  the PRIN 2022 project \textit{High order structure-preserving semi-implicit schemes for hyperbolic equations} and via the  Departments of Excellence  Initiative 2018--2027 attributed to DICAM of the University of Trento (grant L. 232/2016).  
MD received further funding by the Fondazione Caritro via the project SOPHOS and by the European Union Next Generation EU project PNRR Spoke 7 CN HPC 
and by the European Research Council (ERC) under the European Union’s Horizon 
2020 research and innovation programme, Grant agreement No. ERC-ADG-2021-101052956-BEYOND. 

\bibliographystyle{plain}
\bibliography{biblio}
\end{document}